\documentclass[smallextended]{svjour3}
\usepackage{amsmath}
\usepackage{stmaryrd}
\usepackage{bbding}
\usepackage{dcolumn}
\usepackage{graphicx}
\usepackage{amsfonts}
\usepackage{amssymb}
\usepackage{psfrag}
\usepackage{wrapfig}
\usepackage{subfigure}
\usepackage{makeidx}
\usepackage{bm}
\usepackage{epsf}
\usepackage{epsfig}
\usepackage{setspace}
\usepackage{color}

\begin{document}
\title{High-order gas-kinetic scheme for radiation hydrodynamics in equilibrium-diffusion limit}
\author{Yaqing Yang \and Liang Pan  \and Wenjun Sun}
\institute{
Yaqing Yang
\at Laboratory of Mathematics and Complex Systems, School of Mathematical Sciences, Beijing Normal University, Beijing, China \\
\email{yqyangbnu@163.com}
\and
Liang Pan 
\at Laboratory of Mathematics and Complex Systems, School of Mathematical Sciences, Beijing Normal University, Beijing, China \\
\email{panliang@bnu.edu.cn}
\and
Wenjun Sun \at
Center for Applied Physics and Technology, College of Engineering, Peking University, Beijing, China\\
Institute of Applied Physics and Computational Mathematics, Beijing, China\\
\email{sun\_wenjun@iapcm.ac.cn}}

\date{Received: date/Accepted: date}

\maketitle
\begin{abstract}
In this paper, a high-order gas-kinetic scheme is developed for the
equation of radiation hydrodynamics in equilibrium-diffusion limit
which describes the interaction between matter and radiation. To
recover RHE, the Bhatnagar-Gross-Krook (BGK) model with modified
equilibrium state is considered. In the equilibrium-diffusion limit,
the time scales of radiation diffusion and hydrodynamic part are
different, and it will make the time step very small for the fully
explicit scheme. An implicit-explicit (IMEX) scheme is applied, in
which the hydrodynamic part is treated explicitly and the radiation
diffusion is treated implicitly. For the hydrodynamics part, a
time dependent gas distribution function can be constructed by the
integral solution of modified BGK equation, and the time dependent
numerical fluxes can be obtained by taking moments of gas
distribution function. For the radiation diffusion term, the
nonlinear generalized minimal residual (GMRES) method is used. To
achieve the temporal accuracy, a two-stage method is developed,
which is an extension of two-stage method for hyperbolic
conservation law. For the spatial accuracy,  the multidimensional
weighted essential non-oscillation (WENO) scheme is used for the
spatial reconstruction. A variety  of numerical tests are provided
for the performance of current scheme, including the order of
accuracy and robustness.

\keywords{High-order gas-kinetic scheme, equation of radiation hydrodynamics,
GMRES method, WENO scheme.}
\end{abstract}

\maketitle

\section{Introduction}
The equation of radiation hydrodynamics describes the radiative
transport through a fluid with coupled momentum and energy exchange
\cite{Radiation-1,Radiation-2,Radiation-3}. Its applications are
mainly in high-temperature hydrodynamics, including gaseous stars in
astrophysics, supernova explosions, combustion phenomena, reentry
vehicles, fusion physics and inertial confinement fusion. The
importance of thermal radiation increases  as the temperature is
raised in the above problems. Such as for the moderate temperature,
the role of radiation is primarily one of transporting energy by
radiative process. But for the higher temperature, the energy and
momentum densities of the radiation field may become comparable to
or even dominates the corresponding fluid quantities.

In the case of the zero diffusion limit
\cite{radiative-GKS-0,radiative-GKS-1,radiative-GKS-2}, the equation
of radiation hydrodynamics can be written into a nonlinear
hyperbolic system of conservation laws. But for the more complicated
equilibrium-diffusion limit \cite{Radiation-Diffusion-1}, another
nonlinear diffusion term for radiative heat transfer should be
added. Due to the highly non-linearity of this radiation diffusion
terms, it becomes more challenge to design a high-order and robust
numerical method. For the scales of characteristic time between the
radiation and hydrodynamics are different by several orders of
magnitude, and it usually requires the radiation part to be solved
implicitly to guarantee the numerical stability. There are many
numerical method for the radiation hydrodynamics in
equilibrium-diffusion limit. With the operator splitting method, the
Godunov schemes were proposed for the hyperbolic part and an
implicit scheme is proposed for the radiative heat transfer
\cite{Radiation-Diffusion-1,Radiation-Diffusion-2,Radiation-Diffusion-3}.
The only second-order accuracy can be achieved in space and time for
both the equilibrium diffusion and streaming limit, and it is also
capable of computing radiative shock solutions accurately.
Furthermore, to achieve the high-order accuracy,  one-dimension
implicit-explicit (IMEX) Lagrangian high-order scheme was developed
in \cite{Radiation-Diffusion-4}. The essentially non-oscillatory
(ENO) \cite{ENO} method is used for the advection and radiation
diffusion term to obtain the high spatial accuracy, and the strong
stability preserving method is used for high order temporal
accuracy. The more work on the RHEs' computation can be found in
\cite{Radiation-Diffusion-5,Radiation-Diffusion-6,Radiation-Diffusion-7}.

In the last decades, the gas-kinetic scheme (GKS) and based on the
Bhatnagar-Gross-Krook (BGK) model \cite{BGK-1,BGK-2} have been
developed systematically for the computations from low speed flows
to supersonic ones \cite{GKS-Xu1,GKS-Xu2}. The gas-kinetic scheme is
based on an analytical integral solution of the BGK equation, and
gas distribution function at a cell interface provides a multi-scale
evolution process from the kinetic particle transport to the
hydrodynamic wave propagation. With the two-stage fourth-order
method for Lax-Wendroff type flow solvers
\cite{GRP-high-1,GRP-high-2}, the high-order gas-kinetic schemes
were developed \cite{GKS-high-1,GKS-high-2}.  The high-order scheme
not only reduces the complexity of computation, but also improves
the accuracy of the numerical solution. Most importantly, the
robustness is as good as the second-order shock capturing scheme.
Furthermore, with the discretization of particle velocity space, a
unified gas-kinetic scheme (UGKS) has been developed for the flow
study in entire Knudsen number regimes from rarefied to continuum
ones \cite{UGKS-Xu1,UGKS-Xu2,UGKS-Xu3}. Recently, the UGKS is
extended to solve radiative transfer system with both scattering and
absorption/emission effects
\cite{radiative-UGKS-1,radiative-UGKS-2,radiative-UGKS-3}. The
asymptotic preserving (AP) property can be accurately recovered. For
the equation of radiation hydrodynamics, a multi-scale scheme is
developed, in which GKS is used for the compressible inviscid flow
and UGKS is used for the non-equilibrium radiative transfer
\cite{radiative-UGKS-4}. Due to the possible large variation of
fluid opacity in different regions, the transport of photons through
the flow system is simulated by the multi-scale scheme.

In this paper, a high-order gas-kinetic scheme is proposed for the
equation of radiation hydrodynamic in the equilibrium-diffusion
limit. Based on the zeroth-order Chapman-Enskog expansion, the
hydrodynamic part of radiation hydrodynamic equation can be obtained
from the modified BGK equation with modified equilibrium state. The
radiation diffusion term is considered as source term.  Since the
time scales of radiation diffusion and hydrodynamic part are
different and it will make the time step of an explicit scheme very
small.  Thus, an implicit-explicit (IMEX) scheme is developed for
solving the radiation hydrodynamic equation, in which the fluid
advection term is treated explicitly and the radiation diffusion is
treated implicitly. For the hydrodynamic part, the gas-kinetic
solver with the modified equilibrium state is used to solve the
compressible flow equations. The nonlinear Newton-GMRES method
\cite{GMRES-1,GMRES-2}  is used to deal with the radiation
diffusion. To achieve the temporal accuracy, the two-stage
third-order temporal discretization is developed, which is an
extension of two-stage method was developed for Lax-Wendroff type
flow solvers. To achieve the spatial accuracy, the classical
weighted essentially non-oscillatory (WENO) \cite{WENO-JS,WENO-Z}
method reconstruction is used. With the two-stage temporal
discretization and WENO reconstruction, a reliable framework was
provided for equation of radiation hydrodynamics. Various numerical
experiments are carried out to validate the performance of current
scheme.

This paper is organized as follows. In Section 2, the equation of
radiation hydrodynamics and corresponding BGK model are introduced.
The high-order gas-kinetic for RHE is presented in Section 3.
Numerical  examples are included in Section 4 and the last section
is the conclusion.

\section{Equation of radiation hydrodynamics and BGK model}

\subsection{Equation of radiation hydrodynamics}
The equation of radiation hydrodynamics (RHE) describes the motion
of flows under a radiation field. It consists of the Euler equations
coupling with the radiation momentum and energy sources and the
radiation-transport equation, and can be given as follows
\begin{align*}
\displaystyle\frac{\partial \rho}{\partial t}&+\nabla\cdot (\rho\boldsymbol{U})=0,\\
\displaystyle\frac{\partial \rho \boldsymbol{U}}{\partial t}&+\nabla\cdot (\rho\boldsymbol{U}\boldsymbol{U}+p)=-\boldsymbol{S}_{rp},\\
\displaystyle\frac{\partial  E}{\partial t}&+\nabla\cdot ((E+p)\boldsymbol{U})=-S_{re},\\
\displaystyle\frac{1}{c}\frac{\partial I_{\nu}}{\partial t}&+\boldsymbol{\Omega}\cdot\nabla I_{\nu}=Q_{\nu},
\end{align*}
where $\rho$, $\boldsymbol{U}$, $E$ and $p$ are the density,
velocity, total energy and pressure of the matter, respectively.
$\bm{S}_{rp}$ is the radiation momentum source, $S_{re}$ is the
radiation energy source, $c$ is the speed of light, $I_{\nu}$ is
radiation intensity, $\bm{\Omega}$ is the photon direction of flight
and $Q_{\nu}$ is the angle and frequency dependent radiation source
representing the radiation-matter interaction. The
radiation-transport equation is essentially the conservation of the
photon number, which reveals the relationship between the photon
free transport and radiation-matter interaction, i.e. photon
emission, photon absorption and photon scatter.  The source terms
$S_{re}$ and $\boldsymbol{S}_{rp}$ can be written as the zeroth and
first order frequency-integrated angular moments of $Q_{\nu}$,
respectively,
\begin{align*}
&S_{re}\equiv\int_{4\pi}\int_{0}^{\infty}Q_{\nu }\mathrm{d}\Omega\mathrm{d}\nu=\frac{\partial \varepsilon}{\partial t}+\nabla\cdot\bm{\mathcal{F}},\\
&\bm{S}_{rp}\equiv\frac{1}{c}\int_{4\pi}\int_{0}^{\infty}\boldsymbol{\Omega} Q_{\nu }\mathrm{d}\Omega\mathrm{d}\nu=\frac{1}{c^2}\frac{\partial \bm{\mathcal{F}}}{\partial t}+\nabla\cdot\mathcal{P},
\end{align*}
where $\varepsilon$ is radiation energy density,
$\boldsymbol{\mathcal{F}}$ is radiation flux and $\mathcal{P}$ is
radiation-pressure. The equilibrium-diffusion approximation imposes
four basic assumptions to simplify RHE, i.e. the photon
mean-free-path is small compared to the size of the
absorption-dominated system, the matter-radiation system is in
thermal equilibrium, the radiation flux is diffusive, and the
radiation pressure is isotropic. It is also assumed that the
radiative temperature and the fluid temperature are equal and the
gas is radiatively opaque so that the equilibrium diffusion will be
dealt with. With the assumptions above, the radiation energy
density, radiation flux and radiation-pressure can be simplified as
\begin{align*}
\varepsilon&=a_RT^4,\\
\boldsymbol{\mathcal{F}}&=-\kappa\nabla T^4+\frac{4}{3}\boldsymbol{U}a_RT^4,\\
\mathcal{P}&=\frac{1}{3}a_RT^4,
\end{align*}
where $T$ is temperature, $\kappa$ is diffusion constant and $a_R$
is radiation constant representing the ratio of the radiation energy
to the material energy. The time-derivative of radiation flux in the
total-momentum can be dropped in accordance with the diffusion
approximation. Therefore, the equation of radiation hydrodynamics in
the equilibrium-diffusion limit can be written as
\begin{equation}\label{RHE}
\begin{split}
\displaystyle\frac{\partial \rho}{\partial t}&+\nabla\cdot (\rho\textbf{U})=0,\\
\displaystyle\frac{\partial \rho \textbf{U}}{\partial t}&+\nabla\cdot (\rho\textbf{U}\textbf{U}+p^*)=0,\\
\displaystyle\frac{\partial  E^*}{\partial t}&+\nabla\cdot ((E^*+p^*)\textbf{U})=\nabla\cdot (\kappa\nabla T^4),
\end{split}
\end{equation}
where the total energy $E^*$ and the total pressure $p^*$ are given
as
\begin{align*}
E^*=&\frac{1}{2}\rho\textbf{U}^2+\frac{p}{\gamma-1}+a_RT^4,\\
p^*=&p+\frac{1}{3}a_RT^4.
\end{align*}
The polytropic ideal gas is considered, and the equation of state is
given by
\begin{align*}
p=(\gamma-1)\rho e=(\gamma-1)\rho c_vT,
\end{align*}
where $e$ is the specific internal energy, $c_v$ is the heat
capacity at constant volume and $\gamma$ is the specific heat ratio.

\subsection{BGK model}
In this paper, a high-order gas-kinetic scheme will be presented for
two-dimensional flows. To recover the macroscopic equation of
radiation hydrodynamics Eq.\eqref{RHE}, the modified BGK equation
\cite{BGK-1,BGK-2} can be written as
\begin{equation}\label{bgk}
f_t+uf_x+vf_y=\frac{g-f}{\tau},
\end{equation}
where $f$ is the gas distribution function, $g$ is the equilibrium
distribution, $\boldsymbol{u}=(u,v)$ is the particle velocity and
$\tau$ is the collision time. For the equation of radiation
hydrodynamics, a modified equilibrium state function is introduced
\cite{radiative-GKS-2}
\begin{equation*}
\displaystyle
g(\boldsymbol{x},\boldsymbol{u},t)=\rho\big(\frac{\lambda_1\lambda_2}{(\lambda_1+\lambda_2)\pi}\big)^{d/2}
(\frac{\lambda_1}{\pi})^{K_1/2}(\frac{\lambda_2}{\pi})^{K_2/2}
e^{-(\frac{\lambda_1\lambda_2}{\lambda_1+\lambda_2}(\boldsymbol{u}-\boldsymbol{U})^2
+\lambda_1\boldsymbol{\xi}_1^2 +\lambda_2\boldsymbol{\xi}_2^2)},
\end{equation*}
where $d=2$ for two dimensional system, $\boldsymbol{U}$ is the
macroscopic velocity, the internal variables are defined as
$\boldsymbol{\xi}_1^2=(\xi_1)_1^2+…+(\xi_{K_1})_1^2$ and
$\boldsymbol{\xi}_2^2=(\xi_1)_2^2+…+(\xi_{K_2})_2^2$ and the
internal degrees of freedom of $\boldsymbol{\xi}_1$ and
$\boldsymbol{\xi}_2$ satisfy
\begin{equation*}
K_1+d=2/(\gamma-1), ~~K_2+d=6.
\end{equation*}
The parameters $\lambda_1, \lambda_2$ can be given by
\begin{equation*}
\frac{\rho}{2\lambda_1}=p,~~\frac{\rho}{2\lambda_2}=\frac{1}{3} a_RT^4.
\end{equation*}
The collision term also satisfies the compatibility condition
\begin{equation*} 
\int \psi \frac{g-f}{\tau} \text{d}\Xi=0,
\end{equation*}
where
$\displaystyle\psi=(1,u,v,\frac{1}{2}(u^2+v^2+\boldsymbol{\xi}_1^2+\boldsymbol{\xi}_2^2))^T$
and
$\text{d}\Xi=\text{d}u\text{d}v\text{d}\boldsymbol{\xi}_1\text{d}\boldsymbol{\xi}_2$.

According to the Chapman-Enskog expansion, the Euler and
Navier-Stokes equations can be derived form the BGK equation
\cite{GKS-Xu1,GKS-Xu2}. Similarly, the macroscopic equations
Eq.\eqref{RHE} can be also derived from Eq.\eqref{bgk}. Taking
zeroth-order Chapman-Enskog expansion, i.e. with $f=g$, and taking
moments of BGK equation Eq.\eqref{bgk},  we have
\begin{equation*}
\int\psi (g_t+ug_x+vg_y)\text{d}\Xi=0,
\end{equation*}
where
\begin{align*}
Q=\int\psi g\text{d}\Xi=\left(
\begin{array}{c}
\rho\\
\rho U  \\
\rho V \\
\displaystyle\frac{1}{2}\rho(U^2+V^2+\frac{K_1+d}{2\lambda_1}+\frac{K_2+d}{2\lambda_2}) \\
\end{array}\right),
\end{align*}
and
\begin{align*}
F(Q)=\int\psi u g\text{d}\Xi=\left(
\begin{array}{c}
\rho U\\
\displaystyle \rho U^2+\frac{\rho}{2\lambda_1}+\frac{\rho}{2\lambda_2}\\
\rho UV \\
\displaystyle\frac{1}{2}\rho U(U^2+V^2+\frac{K_1+d+2}{2\lambda_1}+\frac{K_2+d+2}{2\lambda_2}) \\
\end{array}\right).
\end{align*}
According to the definition of $\lambda_1, \lambda_2$ and the relation
of $K_1, K_2$, the hyperbolic part of macroscopic equations
Eq.\eqref{RHE} can be recovered and the vector form is used in the
following sections
\begin{equation*}
\displaystyle\frac{\partial Q}{\partial t}+\frac{\partial
F(Q)}{\partial x}+\frac{\partial G(Q)}{\partial y}=S(Q),
\end{equation*}
where $Q$ is the conservative variable, $F(Q)$ and $G(Q)$ are the
fluxes in $x$ and $y$ directions and $S(Q)$ is the source term for
radiation diffusion.

\section{High-order gas-kinetic scheme}
\subsection{Temporal discretization}
Recently, based on the time-dependent flux function of the
generalized Riemann problem solver (GRP) \cite{GRP-high-1,GRP-high-2} and 
gas-kinetic scheme (GKS) \cite{GKS-high-1,GKS-high-2}, a two-stage fourth-order time-accurate discretization
was developed for Lax-Wendroff type flow solvers, particularly
applied for the hyperbolic conservation laws. Considering the following
time-dependent equation with the initial condition
\begin{align*}
\frac{\text{d} Q_{ij}}{\text{d} t}=\mathcal {L}(Q_{ij}),
\end{align*}
where $\mathcal {L}$ is an operator for spatial derivative of flux.
Introducing an intermediate state at $t^*=t_n+\Delta t/2$, the
two-stage temporal discretization can be written as
\begin{align}\label{two-stage}
\begin{split}
&Q_{ij}^*=Q_{ij}^n+\frac{1}{2}\Delta t\mathcal{L}(Q_{ij}^n)+\frac{1}{8}\Delta t^2\frac{\partial}{\partial t}\mathcal{L}(Q_{ij}^n),\\
Q_{ij}^{n+1}=&Q_{ij}^n+\Delta t\mathcal{L}(Q_{ij}^n)+\frac{1}{6}\Delta t^2\big(\frac{\partial}{\partial t}\mathcal{L}(Q_{ij}^n)+2\frac{\partial}{\partial t}\mathcal{L}(Q_{ij}^*)\big).
\end{split}
\end{align}
It can be proved that for the hyperbolic equations the two-stage
time stepping method Eq.\eqref{two-stage} provides a fourth-order
time accurate solution for $Q(t)$ at $t=t_n +\Delta t$. Based on the
high-order spatial reconstruction \cite{WENO-JS,WENO-Z}, successes
have also been achieved for the construction of high-order
gas-kinetic scheme for Euler and Navier-Stokes equations
\cite{GKS-high-1,GKS-high-2}. The two-stage method provides a
reliable framework to develop high-order scheme with the
implementation of second-order flux function. Most importantly, due
to the use of both flux function and its temporal derivative, this
scheme is robust and works perfectly from the subsonic to the
hypersonic  flows.

In this paper, the high-order gas-kinetic scheme will be developed
for the radiation hydrodynamics equation Eq.\eqref{RHE} as well. For
simplicity, the two-dimensional uniform mesh is used. Taking moments
of the BGK equation Eq.\eqref{bgk} and integrating with respect to
space for the cell
$I_{ij}=[x_{i-1/2},x_{i+1/2}]\times[y_{j-1/2},y_{j+1/2}]$, the
semi-discrete finite volume scheme can be written as
\begin{align}\label{source}
\frac{\text{d}Q_{ij}}{\text{d}t}&=\mathcal {L}(Q_{ij})+\mathcal {S}(Q_{ij}),
\end{align}
where $Q_{ij}$ is the cell averaged conservative variables over the
cell $I_{ij}$. The operator $\mathcal {L}(Q_{ij})$ for hydrodynamic
part is given by
\begin{align}\label{hydrodynamic}
\mathcal {L}(Q_{ij})=-\frac{1}{\Delta x}(F_{i+1/2,j}(t)-F_{i-1/2,j}(t))-\frac{1}{\Delta y}(G_{i,j+1/2}(t)-G_{i,j-1/2}(t)),
\end{align}
and the operator $\mathcal {S}(Q_{ij})$  is given by
\begin{align*}
\mathcal {S}(Q_{ij})=\frac{1}{\Delta x\Delta y}S(Q_{ij}),
\end{align*}
where $\Delta x$ and $\Delta y$ are the cell size, $F_{i \pm
1/2,j}(t)$ and $G_{i,j \pm 1/2}(t)$ are the time dependent numerical
fluxes at cell interfaces in $x$ and $y$ directions and $S(Q_{ij})$
is the source for radiation diffusion. Without considering the
source term, the explicit scheme is used for the Euler and
Navier-Stokes equations \cite{GKS-high-1}. However, the time scales
of radiation diffusion and fluid advection are different. It will
make the time step of an explicit scheme very small, and the
explicit scheme can be only used efficiently for the mildly
non-relativistic regime for RHE. To improve the efficiency, the RHE
should be discretized in an implicit-explicit (IMEX) procedure, i.e.
the fluid advection term is treated explicitly and the radiation
component is treated implicitly. Similar with Eq.\eqref{two-stage},
introducing an intermediate state $Q^*$ at $t^*=t_n+\Delta t/2$, the
two-stage temporal discretization for Eq.\eqref{source}  is given as
follows
\begin{equation}\label{two-stage2}
\begin{split}
&Q_{ij}^*=Q_{ij}^n+\frac{\Delta t}{2}\mathcal{L}(Q_{ij}^n)+
\frac{\Delta t^2}{8}\mathcal{L}_t(Q_{ij}^n)+ \frac{\Delta
t}{4}\left(\mathcal{S}(Q_{ij}^n)+\mathcal{S}(Q_{ij}^*)\right),\\
Q_{ij}^{n+1}=Q_{ij}^n+&\Delta t\mathcal{L}(Q_{ij}^n)+
\frac{\Delta t^2}{6} (\mathcal{L}_t(Q_{ij}^n)+2\mathcal{L}_t(Q_{ij}^*))
+\frac{\Delta t}{6}(\mathcal{S}(Q_{ij}^{n})+4\mathcal{S}(Q_{ij}^{*})+\mathcal{S}(Q_{ij}^{n+1})).
\end{split}
\end{equation}
Compared with the original two-stage method Eq.\eqref{two-stage},
the trapezoid integration is used for the source terms.

It can be proved that Eq.\eqref{two-stage2} provides a third-order
accurate approximation for the system with sources Eq.\eqref{source}
at $t=t_n+\Delta t$. Integrating Eq.\eqref{source} on the time
interval $[t^n,t^{n+1}]$, we have
\begin{align*}
Q_{ij}^{n+1}-Q_{ij}^n&=\int_{t_n}^{t_n+\Delta t}(\mathcal{L}+\mathcal{S})(Q_{ij}(t))\text{d}t.
\end{align*}
To prove this proposition above, the following Taylor expansion need
to be satisfied  
\begin{align}\label{exp}
\int_{t_n}^{t_n+\Delta
t}(\mathcal{L}+\mathcal{S})(Q_{ij}(t))\text{d}t=&\Delta
t(\mathcal{L}+\mathcal{S})(Q_{ij}^n)+\frac{\Delta
t^2}{2}(\mathcal{L}+\mathcal{S})_t(Q_{ij}^n) \\ \nonumber
+&\frac{\Delta
t^3}{6}(\mathcal{L}+\mathcal{S})_{tt}(Q_{ij}^n)+\mathcal{O}(\Delta
t^4).
\end{align}
According to Eq.\eqref{source} and Cauchy-Kovalevskaya method, the
temporal derivatives can be given by
\begin{align*}
&\mathcal{L}_t=\mathcal{L}_Q(\mathcal {L}+\mathcal{S}),\\
&\mathcal{S}_t=\mathcal{S}_Q(\mathcal{L}+\mathcal{S}),
\end{align*}
and
\begin{align*}
(\mathcal{L}+\mathcal{S})_{tt}=((\mathcal{L}_{QQ}+\mathcal{S}_{QQ})(\mathcal{L}+\mathcal{S})^2+(\mathcal{L}_Q+\mathcal{S}_Q)^2(\mathcal{L}+\mathcal{S})).
\end{align*}
For the operator $\mathcal{S}$, we have the following expansion up
to the corresponding order
\begin{equation}\label{expansion1}
\begin{split}
\displaystyle
\mathcal{S}(Q_{ij}^*)&=\mathcal{S}(Q_{ij}^n)+\mathcal{S}_{Q}(Q_{ij}^*-Q_{ij}^n)
+\frac{\mathcal{S}_{QQ}}{2}(Q_{ij}^*-Q_{ij}^n)^2+\mathcal{O}(Q_{ij}^*-Q_{ij}^n)^3,\\
\displaystyle \mathcal{S}(Q_{ij}^{n+1})=&\mathcal{S}(Q_{ij}^n)
+\mathcal{S}_Q(Q_{ij}^{n+1}-Q_{ij}^n)+\frac{\mathcal{S}_{QQ}}{2}(Q_{ij}^{n+1}-Q_{ij}^n)^2+\mathcal{O}(Q_{ij}^{n+1}-Q_{ij}^n)^3.
\end{split}
\end{equation}
Substituting Eq.\eqref{expansion1} into Eq.\eqref{two-stage2}, we
have
\begin{align}\label{proof2}
Q_{ij}^*-Q_{ij}^n=&\frac{\Delta t}{2}(\mathcal{L}(Q_{ij}^n)+\mathcal{S}(Q_{ij}^n)+\frac{\Delta t^2}{8}\mathcal{L}_t(Q_{ij}^n)
+\frac{\Delta t}{4}\mathcal{S}_Q(Q_{ij}^*-Q_{ij}^n))+\mathcal{O}(\Delta t^3)\nonumber\\
=&\frac{\Delta t}{2}(\mathcal{L}(Q_{ij}^n)+\mathcal{S}(Q_{ij}^n))
+\frac{\Delta t^2}{8}\big(\mathcal{L}_t(Q_{ij}^n)+\mathcal{S}_t(Q_{ij}^n)\big)+\mathcal{O}(\Delta t^3),
\end{align}
and
\begin{align}\label{proof3}
Q_{ij}^{n+1}-Q_{ij}^n=&\Delta t(\mathcal{L}(Q_{ij}^n)+\mathcal{S}(Q_{ij}^n))
+\frac{\Delta t^2}{6}\big(\mathcal{L}_t(Q_{ij}^n)+2\mathcal{L}_t(Q_{ij}^*))\big)\nonumber\\
+&\frac{\Delta t}{6}\big(\mathcal{S}_Q(Q_{ij}^{n+1}-Q_{ij}^n)+4\mathcal{S}_Q(Q_{ij}^*-Q_{ij}^n)\big)\nonumber\\
+&\frac{\Delta t}{12}\big(\mathcal{S}_{QQ}(Q_{ij}^{n+1}-Q_{ij}^n)^2+4\mathcal{S}_{QQ}(Q_{ij}^*-Q_{ij}^n)^2\big).
\end{align}
For the operator $\mathcal{L}$, we have the following expansion up
to the corresponding order as well
\begin{align}\label{expansion2}
\begin{split}
\displaystyle
\mathcal{L}(Q_{ij}^*)&=\mathcal{L}(Q_{ij}^n)+\mathcal{L}_Q(Q_{ij}^*-Q_{ij}^n)
+\frac{\mathcal{L}_{QQ}}{2}(Q_{ij}^*-Q_{ij}^n)^2+\mathcal{O}(Q_{ij}^*-Q_{ij}^n)^2,\\
\displaystyle
\mathcal{L}_Q(Q_{ij}^*)=&\mathcal{L}_Q(Q_{ij}^n)+\mathcal{L}_{QQ}(Q_{ij}^*-Q_{ij}^n)
+\frac{\mathcal{L}_{QQQ}}{2}(Q_{ij}^*-Q_{ij}^n)^2+\mathcal{O}(Q_{ij}^*-Q_{ij}^n)^2.
\end{split}
\end{align}
Substituting Eq.\eqref{proof2} and Eq.\eqref{expansion2} into
Eq.\eqref{proof3}, it is easy to verify
\begin{align*}
Q_{ij}^{n+1}-Q_{ij}^n&=\Delta t(\mathcal{L}+\mathcal{S})(Q_{ij}^n)+\frac{\Delta t^2}{2}(\mathcal{L}_Q+\mathcal{S}_Q)(\mathcal{L}+\mathcal{S})(Q_{ij}^n)\nonumber\\
&+\frac{\Delta t^3}{6}((\mathcal{L}_{QQ}+\mathcal{S}_{QQ})(\mathcal{L}+\mathcal{S})^2+(\mathcal{L}_Q+\mathcal{S}_Q)^2(\mathcal{L}+\mathcal{S}))(Q_{ij}^n)
+\mathcal{O}(\Delta t^4).
\end{align*}
Therefore, the two-stage method Eq.\eqref{two-stage2} provides a
third-order temporal discretization for radiation hydrodynamic
equations.

\subsection{Discretization for hydrodynamic part}
In the following subsections, the implementation of the hydrodynamic
part and radiative part will be given for Eq.\eqref{two-stage2}. For
the hydrodynamic part,  the numerical flux $F_{i+1/2,j}(t)$ in
$x$-direction can be given by Gaussian quadrature
\begin{align}\label{flux_x}
F_{i+1/2,j}(t)=\frac{1}{\Delta
y}\int_{y_{j-1/2}}^{y_{j+1/2}}F_{i+1/2}(y,t)\text{d}y=\sum_{\ell=1}^2\omega_\ell \int\psi u f(x_{i+1/2},y_{j_\ell},t,u,v,\xi)\text{d}\Xi,
\end{align}
where $(x_{i+1/2},y_{j_\ell})$ is the Gaussian quadrature point and
$\omega_\ell$ are quadrature weights. To construct the gas
distribution function $f(x_{i+1/2},y_{j_\ell},t,u,v,\xi)$ at the
cell interface, the integral solution of BGK equation Eq.\eqref{bgk}
is used
\begin{equation}\label{integral1}
f(x_{i+1/2},y_{j_\ell},t,u,v,\xi)=\frac{1}{\tau}\int_0^t g(x',y',t',u,v,\xi)e^{-(t-t')/\tau}dt'+e^{-t/\tau}f_0(-ut,-vt,u,v,\xi),
\end{equation}
where $f_0$ is the initial gas distribution function, $g$ is the
corresponding equilibrium state, and $x_{i+1/2}=x'+u(t-t')$ and
$y_{j_\ell}=y'+v(t-t')$ are the trajectory of particles. Similar
with the gas-kinetic scheme for Euler and Navier-Stokes equations,
the  second-order gas-kinetic solver \cite{GKS-Xu2} can be written
as follows
\begin{align}\label{2nd-simplify-flux}
f(x_{i+1/2},y_{j_\ell},t,u,v,\xi)=&(1-e^{-t/\tau})g_0+((t+\tau)e^{-t/\tau}-\tau)(\overline{a}u+\overline{b}v)g_0\nonumber\\
+&(t-\tau+\tau e^{-t/\tau}){\bar{A}} g_0\nonumber\\
+&e^{-t/\tau}g_r[1-(\tau+t)(a^{r}u+b^{r}v)-\tau A^r)](1-H(u))\nonumber\\
+&e^{-t/\tau}g_l[1-(\tau+t)(a^{l}u+b^{l}v)-\tau A^l)]H(u).
\end{align}
The coefficients in Eq.\eqref{2nd-simplify-flux} can be determined
by the reconstructed derivatives and compatibility condition
\begin{align*}
\displaystyle
\langle a^{l,r}\rangle=\frac{\partial Q_{l,r}}{\partial x},
\langle b^{l,r}\rangle=\frac{\partial Q_{l,r}}{\partial y},
\langle a^{l,r}u+b^{l,r}v+A^{l,r}\rangle=0,
\end{align*}
and
\begin{align*}
\displaystyle
\langle\overline{a}\rangle=\frac{\partial Q_{0}}{\partial x},
\langle\overline{b}\rangle=\frac{\partial Q_{0}}{\partial y},
\langle\overline{a}u+\overline{b}v+\overline{A}\rangle=0,
\end{align*}
where the moments of the equilibrium $g$ are defined by
\begin{align*}
\langle...\rangle=\int g (...)\psi \text{d}\Xi.
\end{align*}
Compared with the Euler and Navier-Stokes equations, the procedure
for RHE is a little more complicated. As an example,  the spatial
derivative of $g$ can be written as
\begin{equation}\label{equi-g}
\frac{\partial g}{\partial x}=ag,
\end{equation}
where
\begin{equation*}
a=\boldsymbol{a}\cdot\tilde{\boldsymbol{\psi}}=a_1+a_2u+a_3v+\frac{1}{2}a_4(u^2+v^2)+\frac{1}{2}a_5\boldsymbol{\xi}_1^2+\frac{1}{2}a_6\boldsymbol{\xi}^2_2.
\end{equation*}
Taking moments of Eq.\eqref{equi-g}, we have
\begin{equation*}
\frac{\partial Q}{\partial x}=\int \frac{\partial g}{\partial x}\boldsymbol{\psi}\text{d}\Xi=\int ag\boldsymbol{\psi}\text{d}\Xi,
\end{equation*}
To obtain the connections $\boldsymbol{a}=(a_1,…,a_6)$, the
equation can be given as a linear system
\begin{equation}\label{linear}
M\boldsymbol{a}=\frac{\partial Q}{\partial x},
\end{equation}
where
\begin{equation*}
\displaystyle M=\int\boldsymbol{\psi}\otimes\tilde{\boldsymbol{\psi}}g\text{d}\Xi.
\end{equation*}
The detailed formulation of $M$ can be found in Appendix. The system
Eq.\eqref{linear} seems to be under determined, but we can also find
the unique solution.

According to the chain rule, the derivative of $g$ can be also
written as
\begin{equation*}
\frac{\partial g}{\partial x}=\frac{\partial g}{\partial
\rho}\frac{\partial \rho}{\partial x}+\frac{\partial g}{\partial
U}\frac{\partial U}{\partial x}+\frac{\partial g}{\partial
V}\frac{\partial V}{\partial x} +\frac{\partial g}{\partial
\lambda_1}\frac{\partial \lambda_1}{\partial x}+\frac{\partial
g}{\partial \lambda_2}\frac{\partial \lambda_2}{\partial x}.
\end{equation*}
It's easily to verify that
\begin{equation*}
\frac{\partial g}{\partial \lambda_1}= C_1g-\boldsymbol{\xi}_1^2g,~
\frac{\partial g}{\partial \lambda_2}= C_2g-\boldsymbol{\xi}_2^2g,
\end{equation*}
where $C_1$ and $C_2$ are the functions without
$\boldsymbol{\xi}_1^2$ and $\boldsymbol{\xi}_2^2$, respectively.
Comparing the coefficients of $\boldsymbol{\xi}_1^2$ and
$\boldsymbol{\xi}_2^2$ with Eq.\eqref{equi-g}, we have
\begin{equation*}
\frac{1}{2}a_5\boldsymbol{\xi}_1^2g=-\boldsymbol{\xi}_1^2g\frac{\partial
\lambda_1}{\partial x},~
\frac{1}{2}a_6\boldsymbol{\xi}_2^2g=-\boldsymbol{\xi}_2^2g\frac{\partial
\lambda_2}{\partial x}.
\end{equation*}
According to the definition of $ \lambda_1$ and $\lambda_2$, we have
\begin{align*}
a_5&=-2\frac{\partial \lambda_1}{\partial x}=\frac{T_x}{RT^2},\\
a_6&=-2\frac{\partial \lambda_2}{\partial x}=-\frac{\rho_xT-4\rho T_x}{a_RT^5/3},
\end{align*}
where
\begin{align*}
\displaystyle T_x&=\frac{1}{\rho c_v+4a_RT^3}((E^*_x-\frac{(\rho
U)(\rho U)_x+(\rho V)(\rho V)_x}{\rho})+\frac{((\rho U)^2+(\rho
V)^2)\rho_x}{2\rho^2}-\rho_xc_vT).
\end{align*}
With the solution for $a_5$ and $a_6$, the unique solution of the
linear system is given by
\begin{align*}
a_4&=\big(\frac{1}{\rho} (A-UB-VC)-\frac{K_1}{8\lambda_1^2}a_5-\frac{K_2}{8\lambda_2^2}a_6\big)/L^2,\\
a_3&=\frac{1}{\rho L}C-Va_4,\\
a_2&=\frac{1}{\rho L}B-Ua_4,\\
a_1&=\frac{1}{\rho}\frac{\partial\rho}{\partial x}-Ua_2-Va_3-B_1'a_4-\frac{K_1}{4 \lambda_1}a_5-\frac{K_2}{4 \lambda_2}a_6,
\end{align*}
where
\begin{align*}
A&=\frac{\partial E^*}{\partial x}-\frac{1}{2}(U^2+V^2+\frac{K_1+2}{2\lambda_1}+\frac{K_2+2}{2\lambda_2})\frac{\partial\rho}{\partial x},\\
B&=\frac{\partial( \rho U)}{\partial x}-U\frac{\partial\rho}{\partial x},\\
C&=\frac{\partial( \rho V)}{\partial x}-V\frac{\partial\rho}{\partial x},\\
L&=\frac{1}{2\lambda_1}+\frac{1}{2\lambda_2}.
\end{align*}
Similarly, the coefficients for temporal derivative can be also
determined. Thus, the gas distribution function
Eq.\eqref{2nd-simplify-flux} and the numerical fluxes
$F_{i+1/2,j}(Q^n,t)$ Eq.\eqref{flux_x} can be fully constructed.

In order to utilize the two-stage temporal discretization, the
temporal derivatives of the flux function need to be determined. The
flux $F_{i+1/2,j}(Q^n,t)$ in the time interval $[t_n, t_n+\Delta t]$
is expanded as the following linear form
\begin{align}\label{expansion}
F_{i+1/2,j}(t)=F_{i+1/2,j}^n+ \partial_t F_{i+1/2,j}^n(t-t_n).
\end{align}
where the coefficients $F_{i+1/2,j}^n$ and $\partial_tF_{i+1/2,j}^n$
can be fully determined by solving the linear system
\begin{align*}
F_{i+1/2,j}^n\Delta t&+\frac{1}{2}\partial_t
F_{i+1/2,j}^n\Delta t^2 =\int_{t_n}^{t_n+\Delta t}F_{i+1/2,j}(t)dt, \\
\frac{1}{2}F_{i+1/2,j}^n\Delta t&+\frac{1}{8}\partial_t
F_{i+1/2,j}^n\Delta t^2 =\int_{t_n}^{t_n+\Delta t/2}F_{i+1/2,j}(t)dt.
\end{align*}
Similarly, the coefficients $G_{i,j+1/2}^n$ and
$\partial_tG_{i,j+1/2}^n$ corresponding to the flux $G_{i,j+1/2}(t)$
in $y$-direction can be constructed as well. According to
Eq.\eqref{hydrodynamic}, $\mathcal {L}(Q_{ij}^n)$ and its temporal
derivative $\partial_t\mathcal {L}(Q_{ij}^n)$ at $t^n$ can be given
by
\begin{align*}
\mathcal{L}(Q_{ij}^n)&=-\frac{1}{\Delta x}
(F_{i+1/2,j}^n(t)-F_{i-1/2,j}^n(t))-\frac{1}{\Delta y}
(G_{i,j+1/2}^n(t)-G_{i,j-1/2}^n(t)),\\
\partial_t\mathcal{L}(Q_{ij}^n)=&-\frac{1}{\Delta x}
(\partial_tF_{i+1/2,j}^n(t)-\partial_tF_{i-1/2,j}^n(t))-\frac{1}{\Delta y}
(\partial_tG_{i,j+1/2}^n(t)-\partial_tG_{i,j-1/2}^n(t)).
\end{align*}
With the procedure at the intermediate state, $\mathcal
{L}(Q_{ij}^*), \partial_t\mathcal{L}(Q_{ij}^*)$  can be constructed
as well.

\subsection{Discretization for radiative part}
According to the definition of source terms, we only need to solve
energy equations at two-stages for Eq.\eqref{two-stage2}, which are
nonlinear systems with respect to the temperature $T$. At each
stage, the energy equations  can be simplified as the nonlinear
system
\begin{equation}\label{orinonlinear}
\mathcal{F}(T)=0,
\end{equation}
where $\mathcal{F}$ is a nonlinear function from $R^N$ to $R^N$ and
$N$ is the number of cells. To obtain an approximate solution of the
nonlinear system, the above system is rewritten as
\begin{equation*}
\mathcal{J}\boldsymbol{\delta}=-\mathcal{F},
\end{equation*}
where $\mathcal{J}$ is the Jacobian of $\mathcal{F}$. In order to
avoid calculating the Jabobian directly, $\mathcal{J}$ can be
approximated by the $F$-derivative as follows
\begin{equation}\label{nonlinear}
\displaystyle\mathcal{J}(T)\cdot\boldsymbol{\delta}=\frac{\mathcal{F}(T+ \sigma \boldsymbol{\delta})-\mathcal{F}(T)}{\sigma},
\end{equation}
where $\sigma$ is a small scalar. In fact, it is not trivial to
choose $\sigma$, which has great influence on the stability of
algorithm and improper values will cause blowing up. According to
\cite{GMRES-2}, a special $\sigma$ is given as follows
\begin{equation*}
\sigma=\frac{\sqrt{\eta}}{\parallel\boldsymbol{\delta}\parallel_2^2}\max\{|T\cdot
\boldsymbol{\delta}|,~ \text{typ}T\cdot
|\boldsymbol{\delta}|\}\cdot\mathrm{sign}(T\cdot
\boldsymbol{\delta}),
\end{equation*}
where $\eta$ is the machine epsilon,
$|\boldsymbol{\delta}|=(|\boldsymbol{\delta}_1|,\dots,|\boldsymbol{\delta}_N|)^T$,
$\text{typ}T=(\text{typ}T_1,\dots,\text{typ}T_N)^T$, $\text{typ}T_i$
is the typical size of $T_i$ and the typical size of real number
equals to its order of magnitude plus one. With the initial guess
$T^{(0)}$ for Eq.\eqref{orinonlinear} and
$\boldsymbol{\delta}^{(0)}$ for Eq.\eqref{nonlinear}, the initial
residual is
\begin{equation*}
r^{(0)}=-\mathcal{F}-\mathcal{J}\boldsymbol{\delta}^{(0)},
\end{equation*}
and its generated Krylov subspace is
\begin{equation*}\scriptsize
K_m\equiv
\mathrm{span}\{\boldsymbol{r}^{(0)},\mathcal{J}\boldsymbol{r}^{(0)},\dots,\mathcal{J}^{m-1}\boldsymbol{r}^{(0)}\}.
\end{equation*}
According to the nonlinear Newton-GMRES method \cite{GMRES-2}, an
approximate solution can be given. Besides, in order to improve
computation efficiency, a restart algorithm is applied in
orthogonalization process \cite{GMRES-1}. In actual calculation, the
restart times can be no more than 5 times if the restart step and
convergence condition are set appropriately.

\begin{table}[!h]
\begin{center}
\def\temptablewidth{0.7\textwidth}
{\rule{\temptablewidth}{1.0pt}}
\begin{tabular*}{\temptablewidth}{@{\extracolsep{\fill}}c|cc|cc}
num of cells       & $L^1$ error  &    Order    &  $L^2$ error &  Order    \\
\hline
$10$     &   1.9001E-03  & ~~       &  1.5274E-03  &   ~~      \\
$20$     &   6.3660E-05  & 4.8995   &  5.0075E-05  &  4.9307 \\
$40$     &   2.0139E-06  & 4.9822   &  1.5793E-06  &  4.9867 \\
$80$     &   6.3035E-08  & 4.9977   &  4.9421E-08  &  4.9980 \\
$160$    &   1.9698E-09  & 5.0000   &  1.5443E-09  &  5.0001 \\
\end{tabular*}
{\rule{\temptablewidth}{1.0pt}}
\end{center}
\caption{\label{tab-1d-1} Accuracy test: errors and orders of
accuracy with $\kappa=0$ and $a_R= 0$ at
$t=2$.}
\begin{center}
\def\temptablewidth{0.5\textwidth}
{\rule{\temptablewidth}{1.0pt}}
\begin{tabular*}{\temptablewidth}{@{\extracolsep{\fill}}c|c|c}
num of cells    & $L^1$ error     &  $L^2$ error    \\
\hline
20     &   1.8422E-03      &  1.4787E-03       \\
40     &   6.1674E-05     &  4.8513E-05   \\
80     &   1.9536E-06     &  1.5314E-06  \\
160    &   6.1704E-08    &  4.8354E-08   \\
\end{tabular*}
{\rule{\temptablewidth}{1.0pt}}
\end{center}
\caption{\label{tab-1d-2} Accuracy test: errors of $\|U_N-U_{N/2}\|$
with $\kappa=0$ and $a_R=10^{-5}$  at $t=2$.}
\begin{center}
\def\temptablewidth{0.5\textwidth}
{\rule{\temptablewidth}{1.0pt}}
\begin{tabular*}{\temptablewidth}{@{\extracolsep{\fill}}c|c|c}
num of cells     & $L^1$ error    &  $L^2$ error    \\
\hline
20      &      1.8407E-03   &  1.4782E-03    \\
40      &      6.1984E-05     &  4.8438E-05    \\
80      &      1.8685E-06     &  1.4935E-06   \\
160    &      8.5963E-08     &  7.5867E-08   \\
\end{tabular*}
{\rule{\temptablewidth}{1.0pt}}
\end{center}
\caption{\label{tab-1d-3} Accuracy test: errors of $\|U_N-U_{N/2}\|$
with $\kappa=10^{-5}$ and $a_R=10^{-5}$ at
$t=2$.}
\end{table}

\section{Numerical tests}
In this section, the one and two dimensional tests cases are
provided to the correctness and robustness of current scheme. In the
computation, the collision time $\tau$ takes
\begin{align*}
\tau=\epsilon \Delta t+C\displaystyle|\frac{p_l-p_r}{p_l+p_r}|\Delta t,
\end{align*}
where $p_l$ and $p_r$ denote the pressure on the left and right
sides of the cell interface. In smooth flow regions, the collision
term takes $\tau=0$ and the gas distribution function
Eq.\eqref{2nd-simplify-flux} will reduce to
\begin{align*}
f=g_0(1+\bar{A}t).
\end{align*}

To achieve the spatial accuracy, the classical multidimensional
fifth-order WENO reconstruction is adopted and more details can be
found in \cite{GKS-high-1,WENO-Z}. In order to eliminate the
spurious oscillation and improve the stability, the WENO
reconstruction is performed for the characteristic variables, and
the detailed analysis for characteristics can be found in
\cite{Radiation-Diffusion-4}. Without special statement, the gas
with $\gamma=5/3$ is used in the following numerical examples.

\begin{figure}[!h]
\centering
\includegraphics[width=0.6\textwidth]{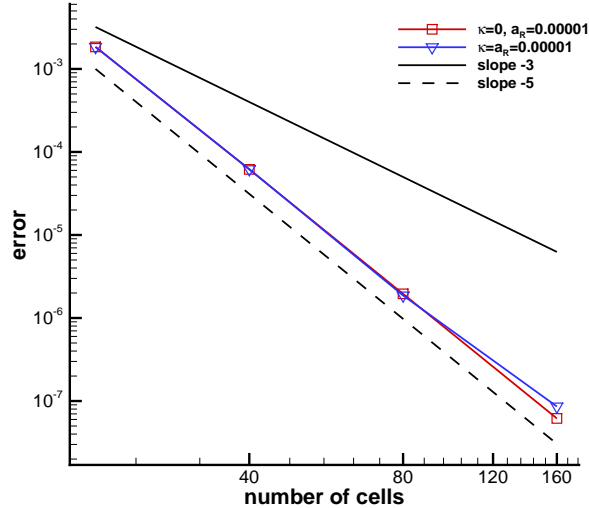}
\caption{\label{accuracy-a} Accuracy test: log-log plots for $L^1$
and $L^2$ norms of $\|U_N-U_{N/2}\|$ and number of cells $N$ at
$t=2$.}
\end{figure}

\subsection{Accuracy test}
The advection of density perturbation problem for Euler equations is
extended to test the order of accuracy for RHE. In this test case,
the initial condition is set as follows
\begin{equation*}
\rho(x)=1+0.2\sin(\pi x),~u(x)=1,~p(x)=1,
\end{equation*}
where  the computational domain is $[0,2]$, the periodic boundary
condition is applied and $c_v=1/(\gamma-1)$.

\begin{figure}[!h]
\centering
\includegraphics[width=0.485\textwidth]{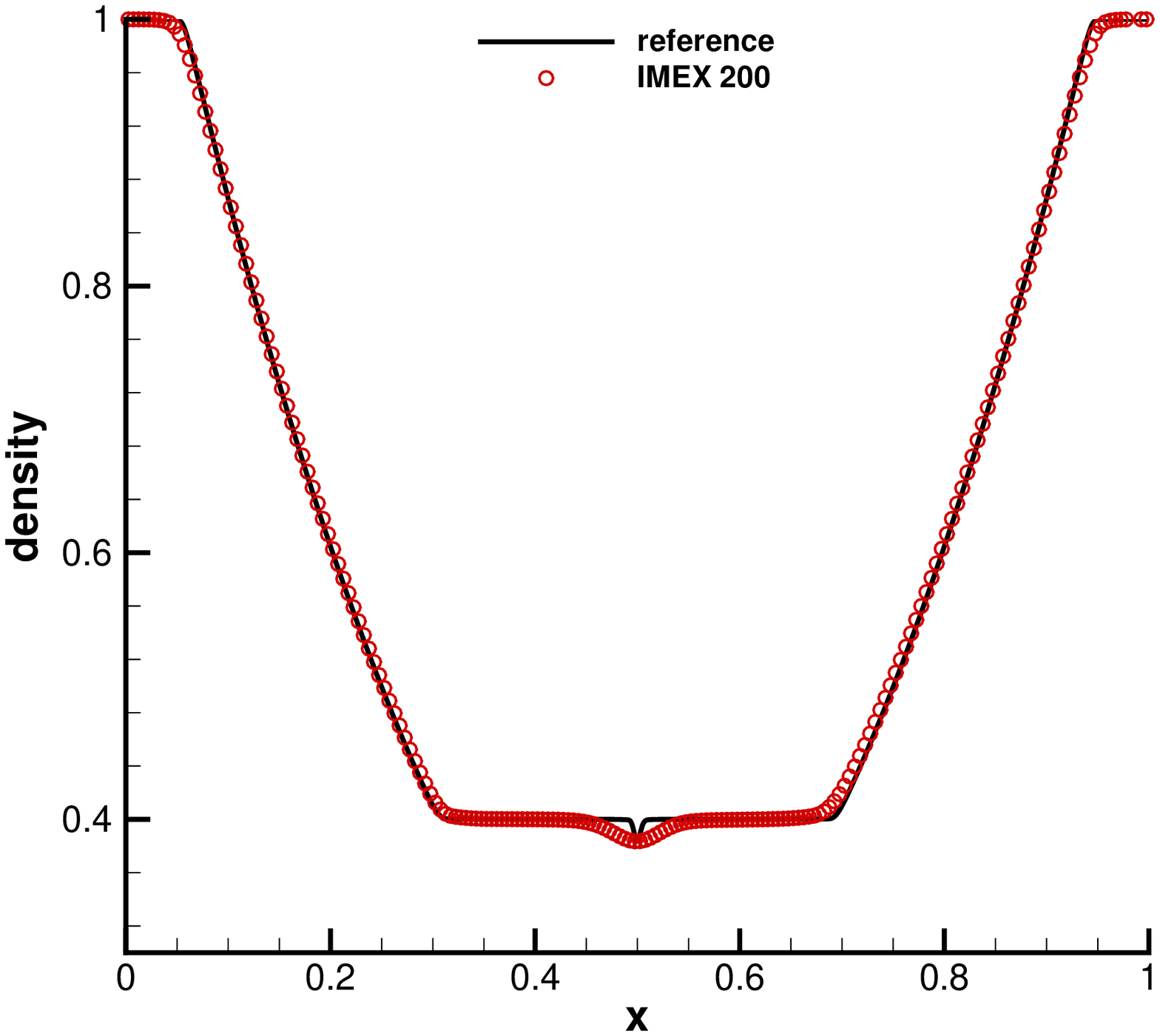}
\includegraphics[width=0.485\textwidth]{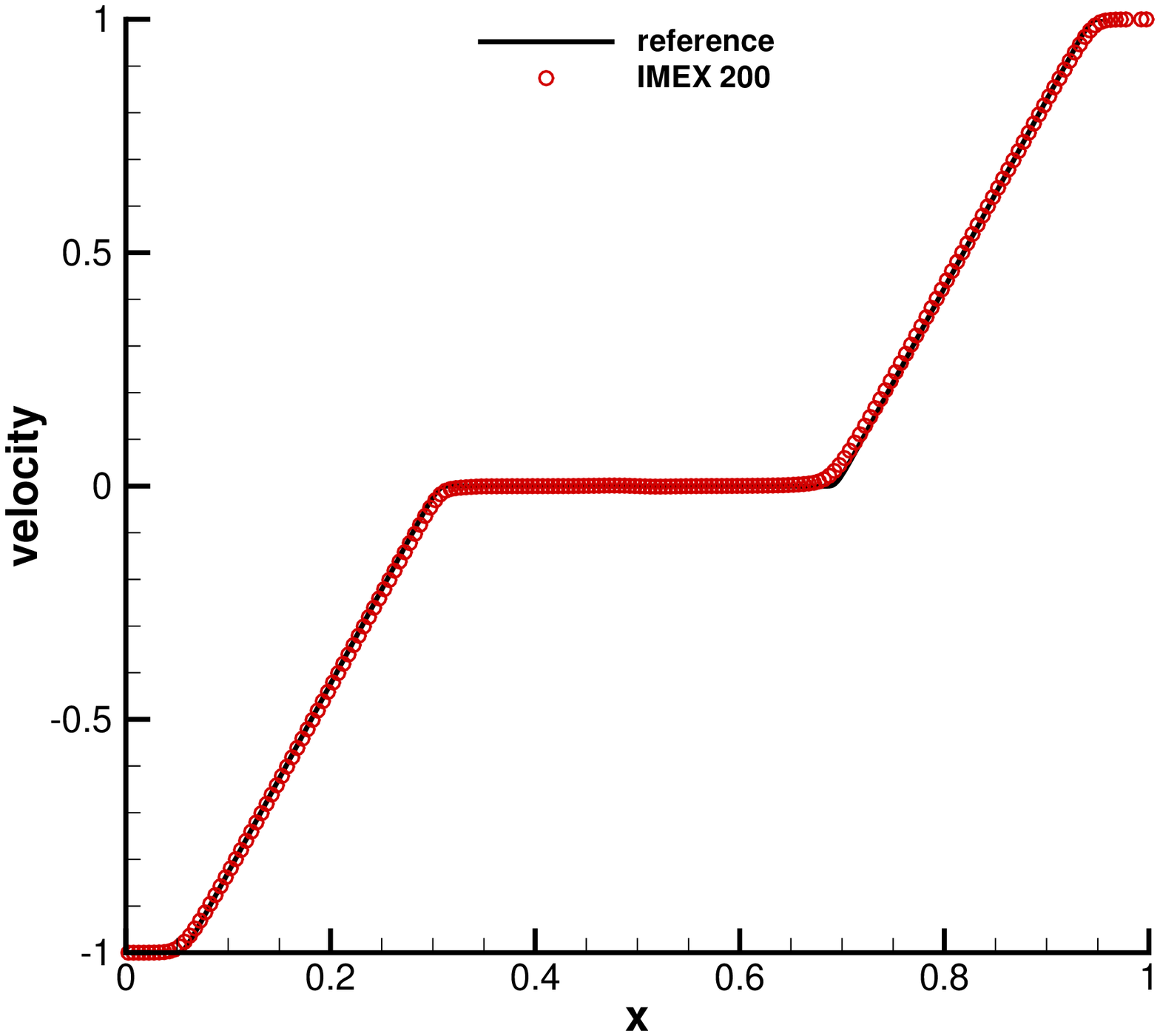}\\
\includegraphics[width=0.485\textwidth]{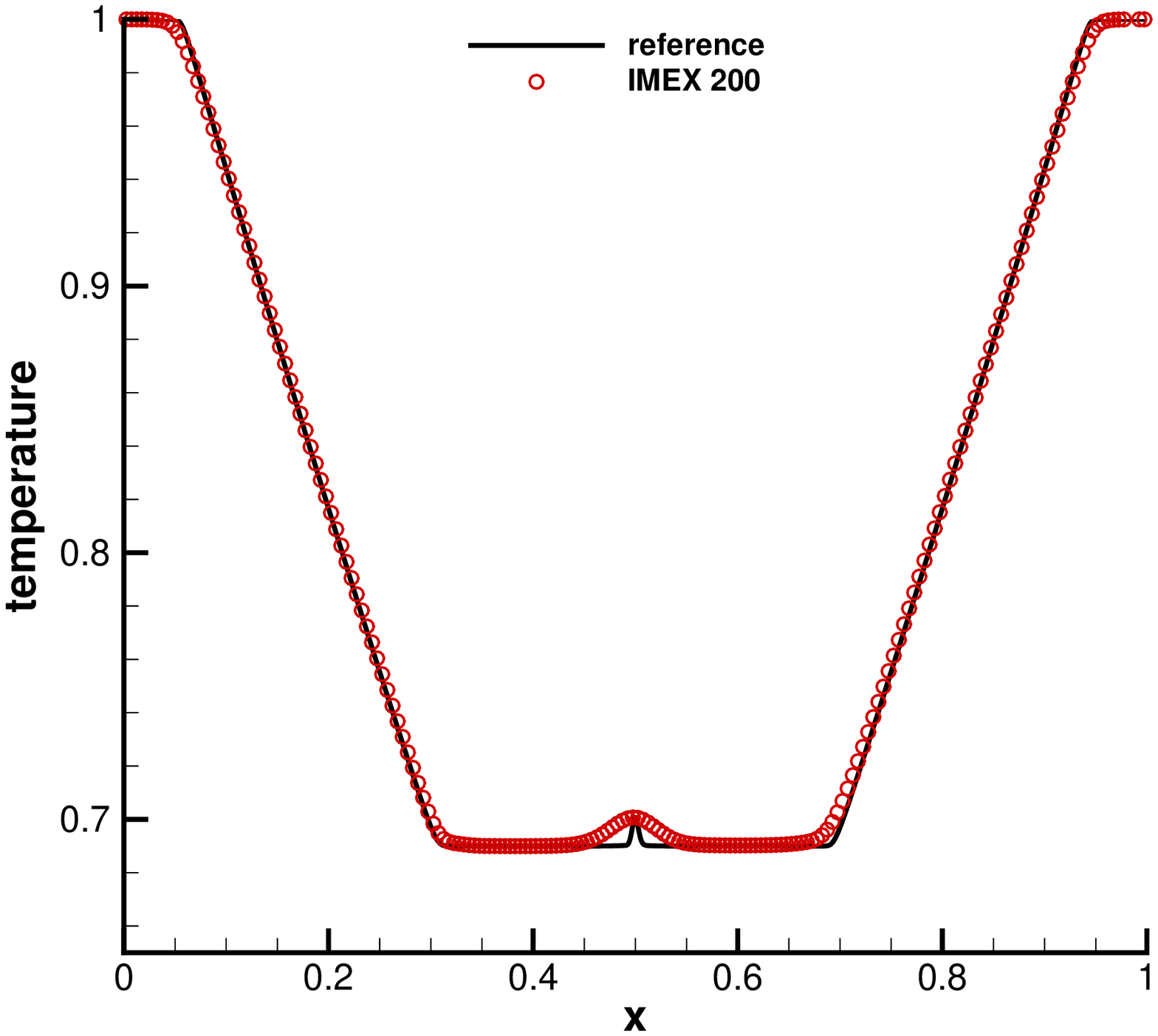}
\includegraphics[width=0.485\textwidth]{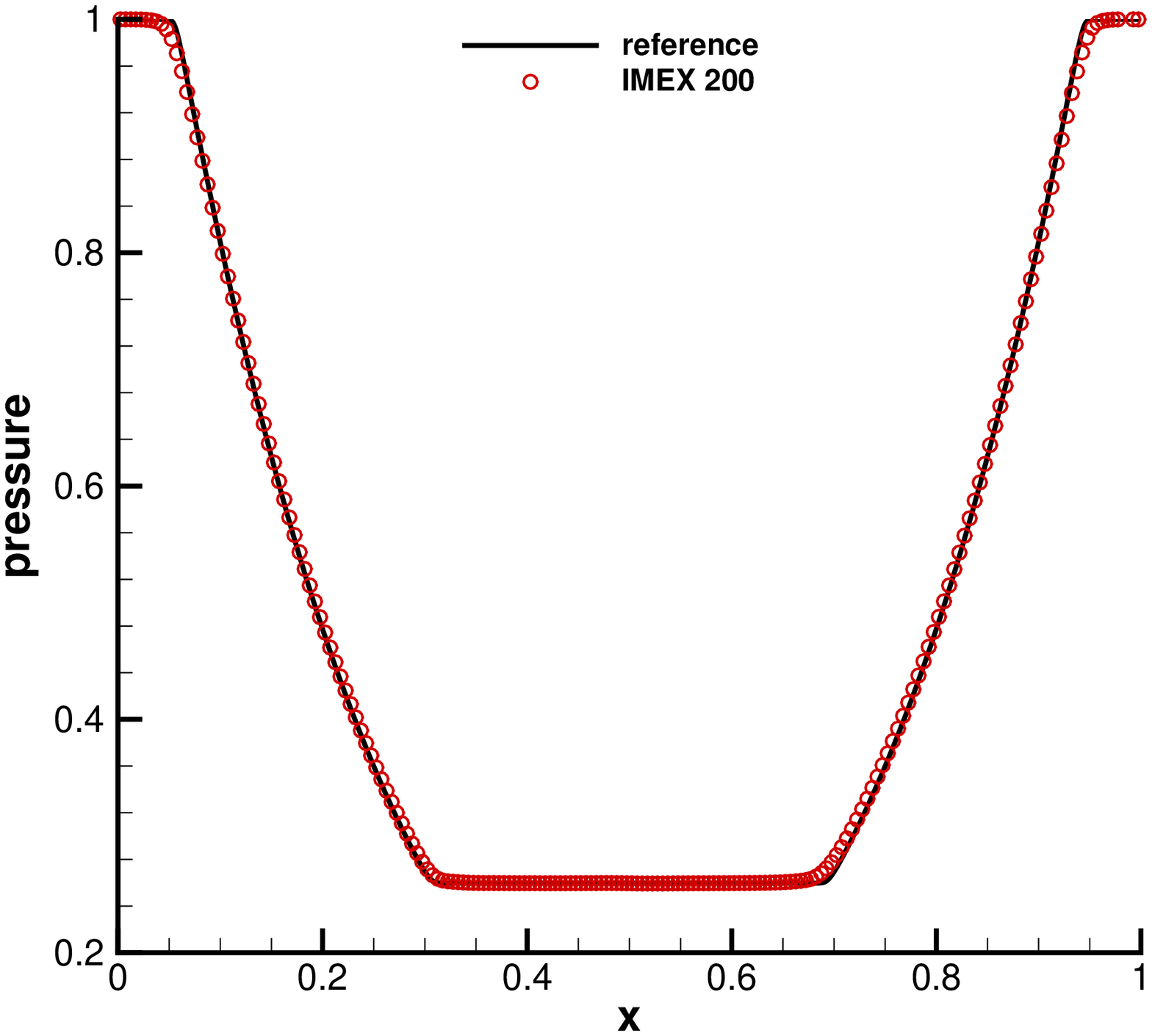}
\caption{\label{riemann-rarefraction-a} Shock tube problem involving
two rarefactions for density, velocity, temperature and pressure at
$t=0.2$ with $\kappa=0$ for the numerical results and reference
data. }
\end{figure}

As reference, the case with $\kappa=0$ and $a_R=0$ is tested. With
these parameters, the RHE degenerates into the Euler equation and
the exact solution can be given by
\begin{equation*}
\rho(x,t)=1+0.2\sin(\pi (x-t)), ~u(x,t)=1, ~p(x,t)=1.
\end{equation*}
The $L^1$ and $L^2$ errors and orders at $t=2$ are shown in
Table.\ref{tab-1d-1} with uniform meshes, where $\Delta x=2/N$ and
$N$ is the number of cells.  The expected order of accuracy can be
achieved. To test the accuracy with radiative effect, two cases with
$\kappa=0, a_R=10^{-5}$ and $\kappa=10^{-5}, a_R=10^{-5}$ are
tested. In the computation, the periodic boundary condition is also
imposed, and the uniform meshes with $N$ cells are used. For these
two cases, they have no exact solution, and  $\|U_N-U_{N/2}\|$ is
computed, where the nonlinear source term is discretized by a
high-order central difference.  For the case with $\kappa=0$, the
system is hyperbolic and the GMRES process is not needed. It falls
into the classical two-stage fourth-order method
Eq.\eqref{two-stage}. For the case with non-zero $\kappa$, the
GMRES process is needed and the tolerance for convergence is set as
$10^{-9}$. The $L^1$ and $L^2$ norms of $\|U_N-U_{N/2}\|$ at $t=2$
are given in Table.\ref{tab-1d-2} and Table.\ref{tab-1d-3}. The
log-log plots for $L^1$ and $L^2$ norms of $\|U_N-U_{N/2}\|$ and
number of cells $N$ is given in Fig.\ref{accuracy-a}. With the
influence of radiation effect and nonlinear source, the theoretical
oder of accuracy can be well kept with the mesh refinement.

\begin{figure}[!h]
\centering
\includegraphics[width=0.485\textwidth]{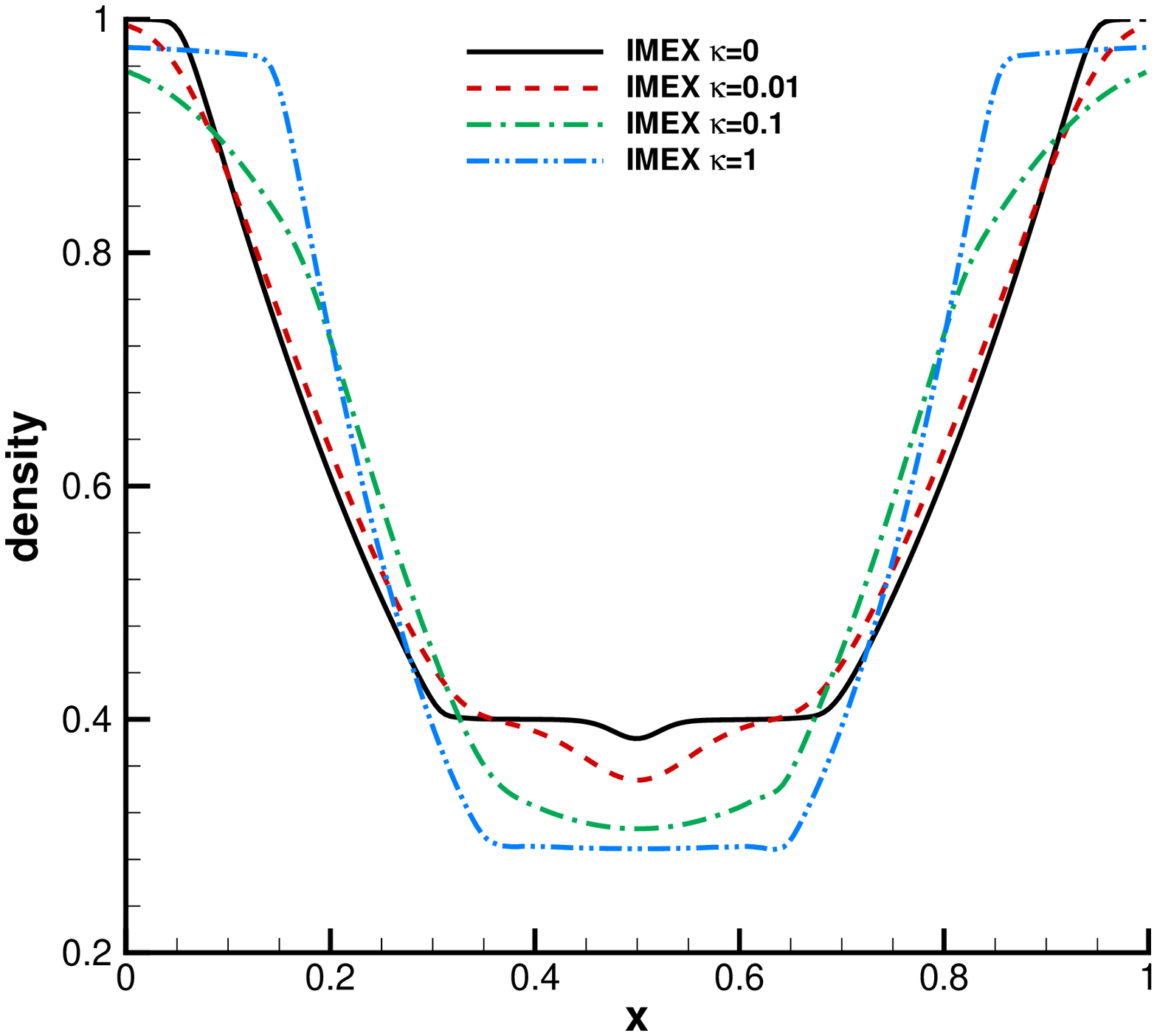}
\includegraphics[width=0.485\textwidth]{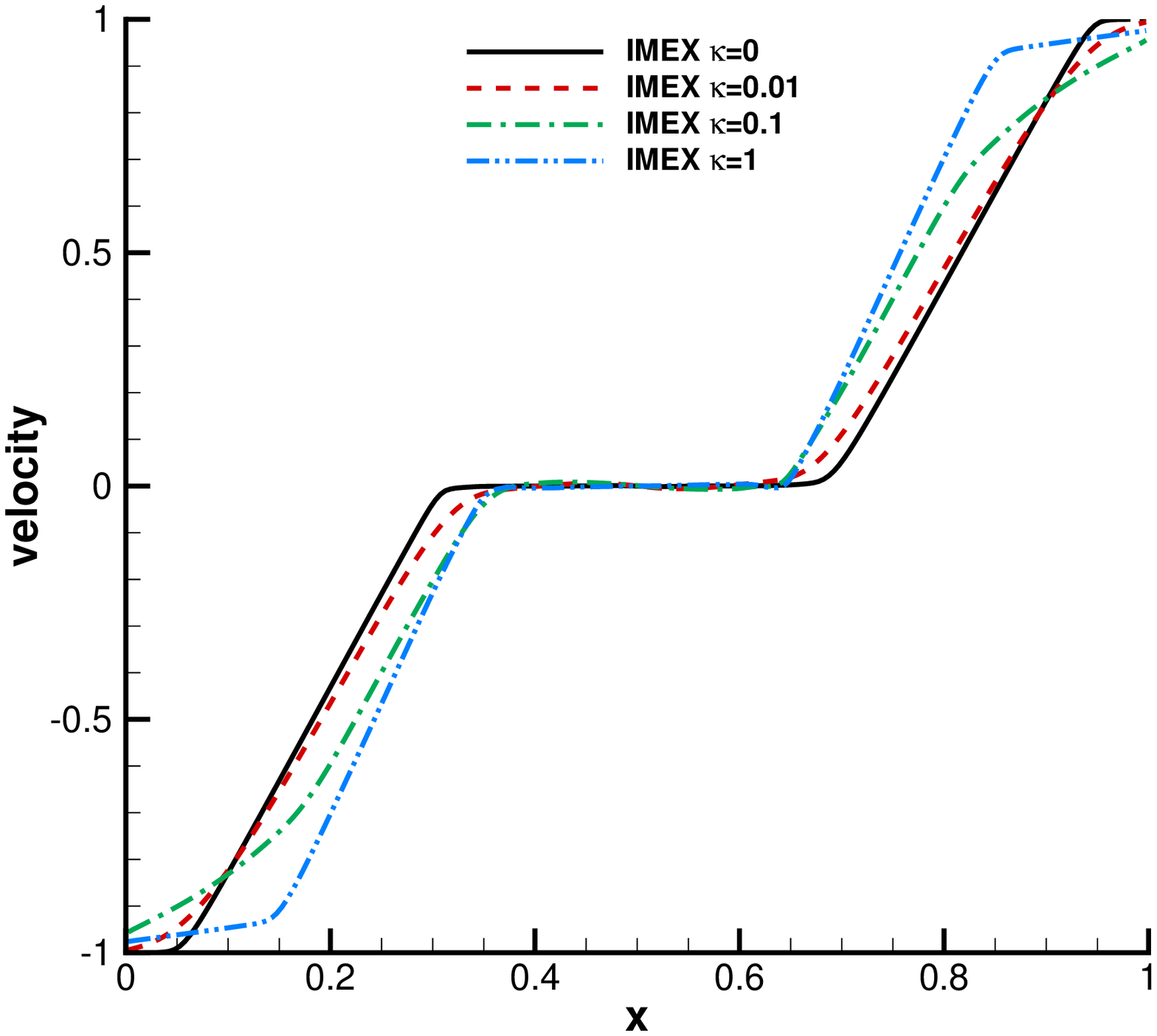}\\
\includegraphics[width=0.485\textwidth]{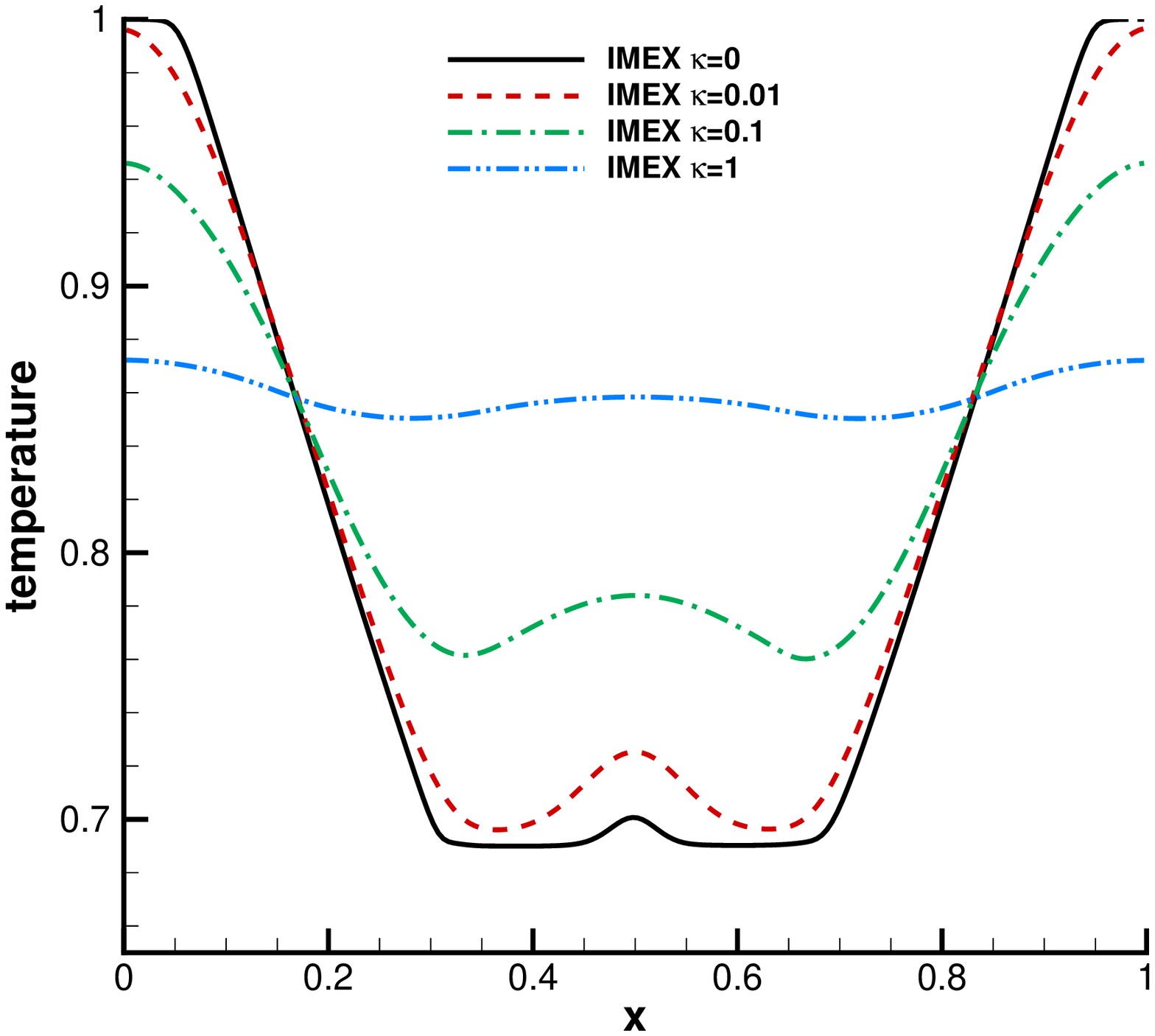}
\includegraphics[width=0.485\textwidth]{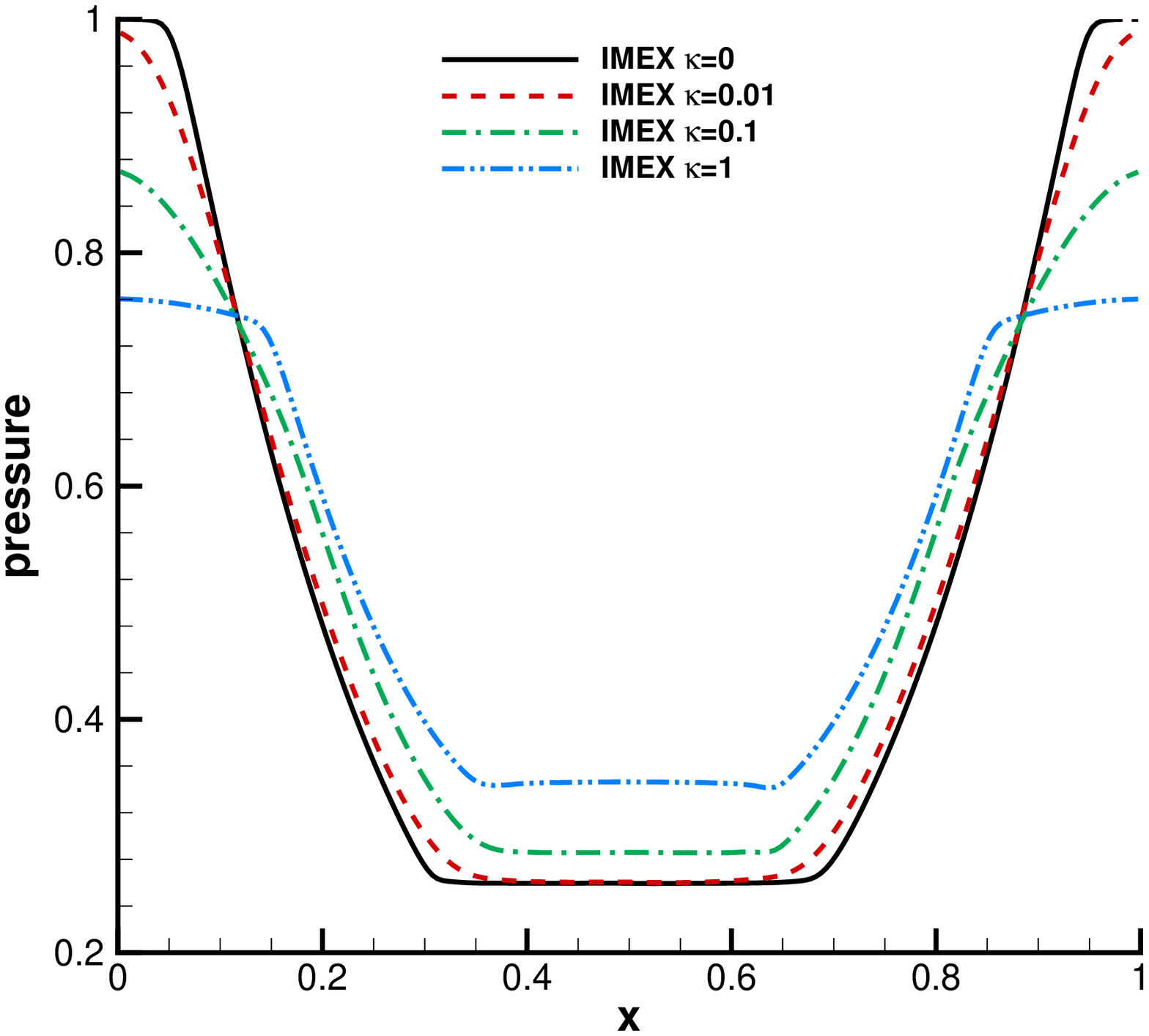}
\caption{\label{riemann-rarefraction-b} Shock tube problem involving two rarefactions:
density, velocity, temperature and pressure at $t=0.2$ with $\kappa=0, 0.01,0.1$ and $1$. }
\end{figure}

\begin{figure}[!h]
\centering
\includegraphics[width=0.485\textwidth]{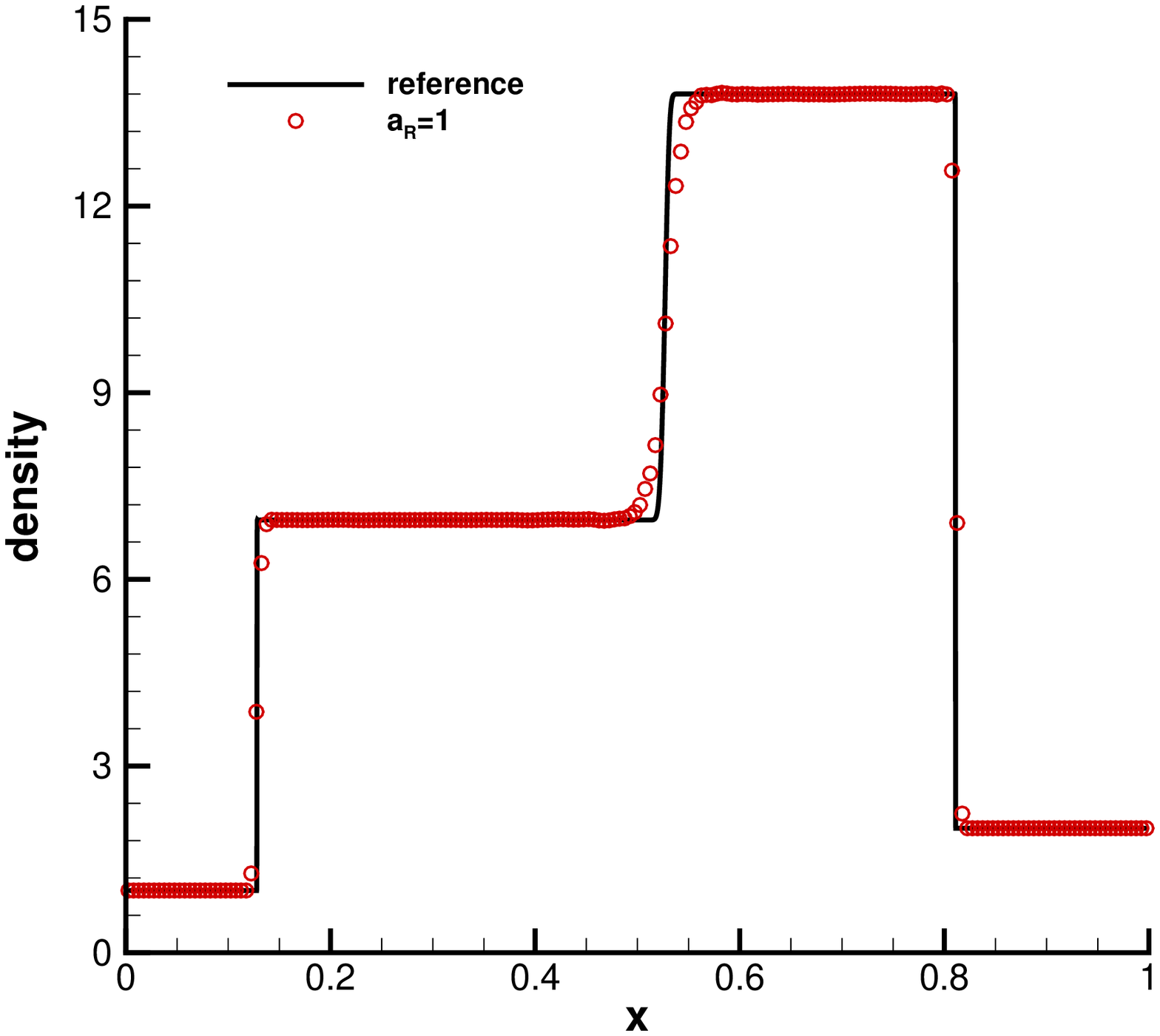}
\includegraphics[width=0.485\textwidth]{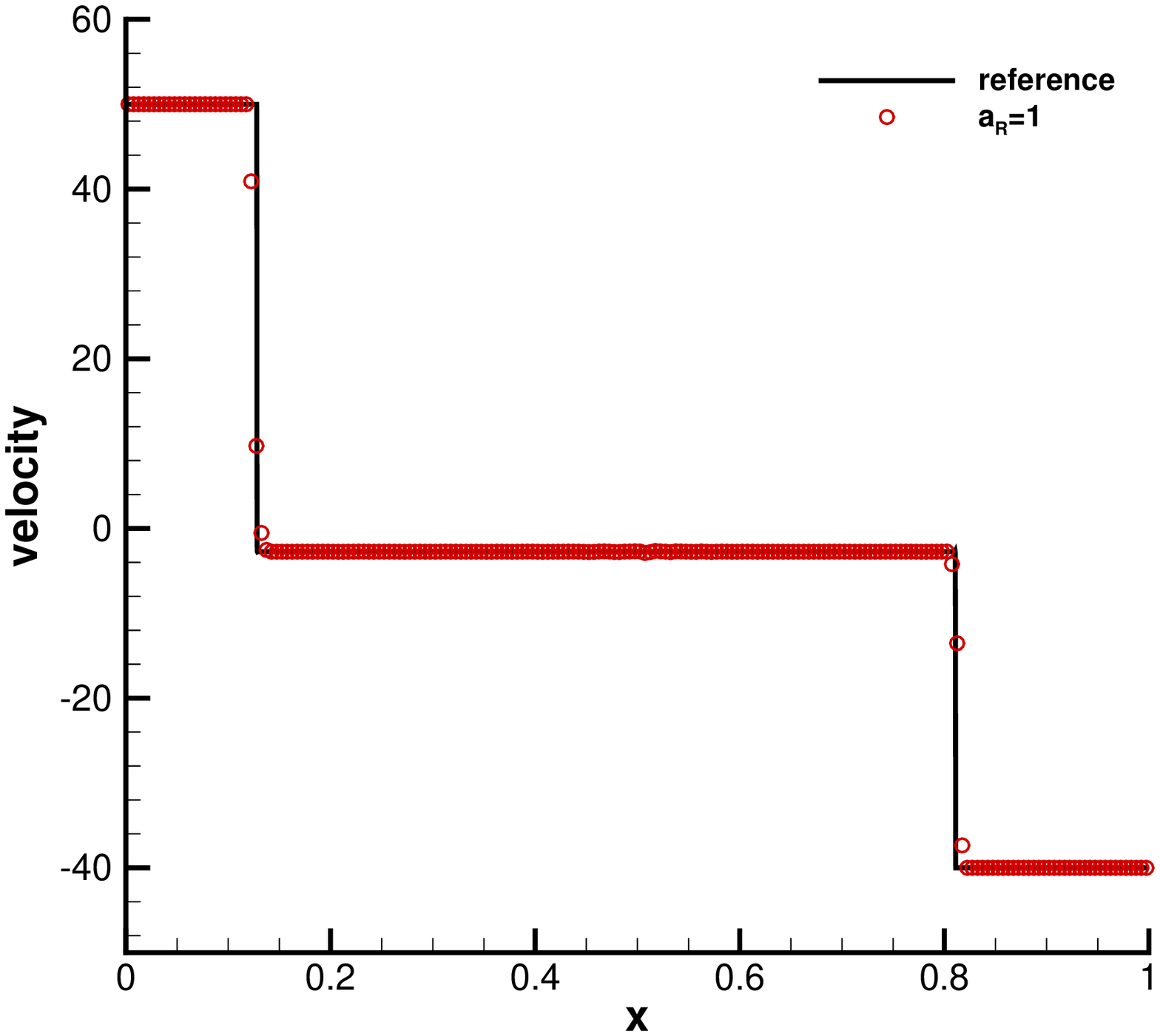}\\
\includegraphics[width=0.485\textwidth]{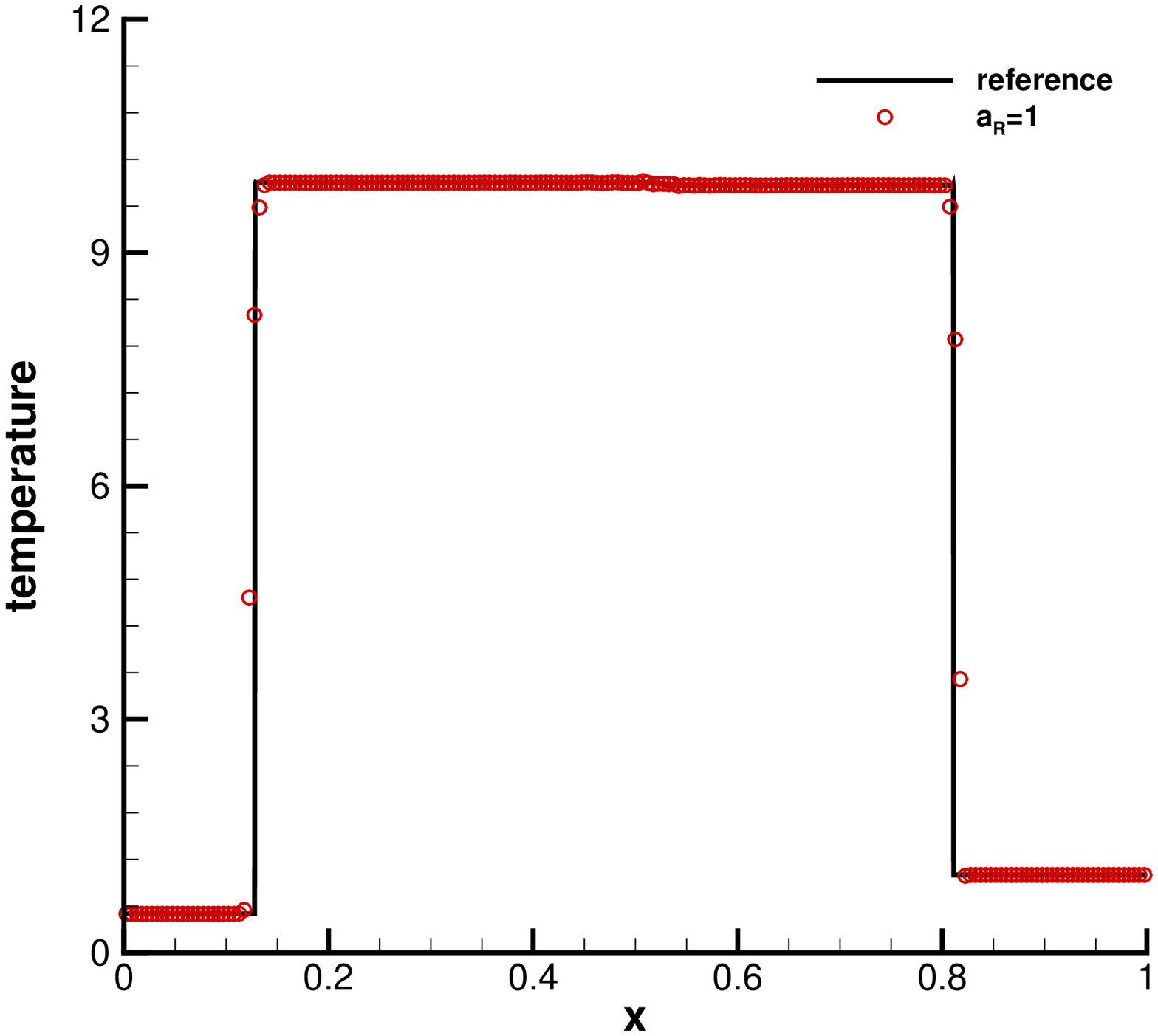}
\includegraphics[width=0.485\textwidth]{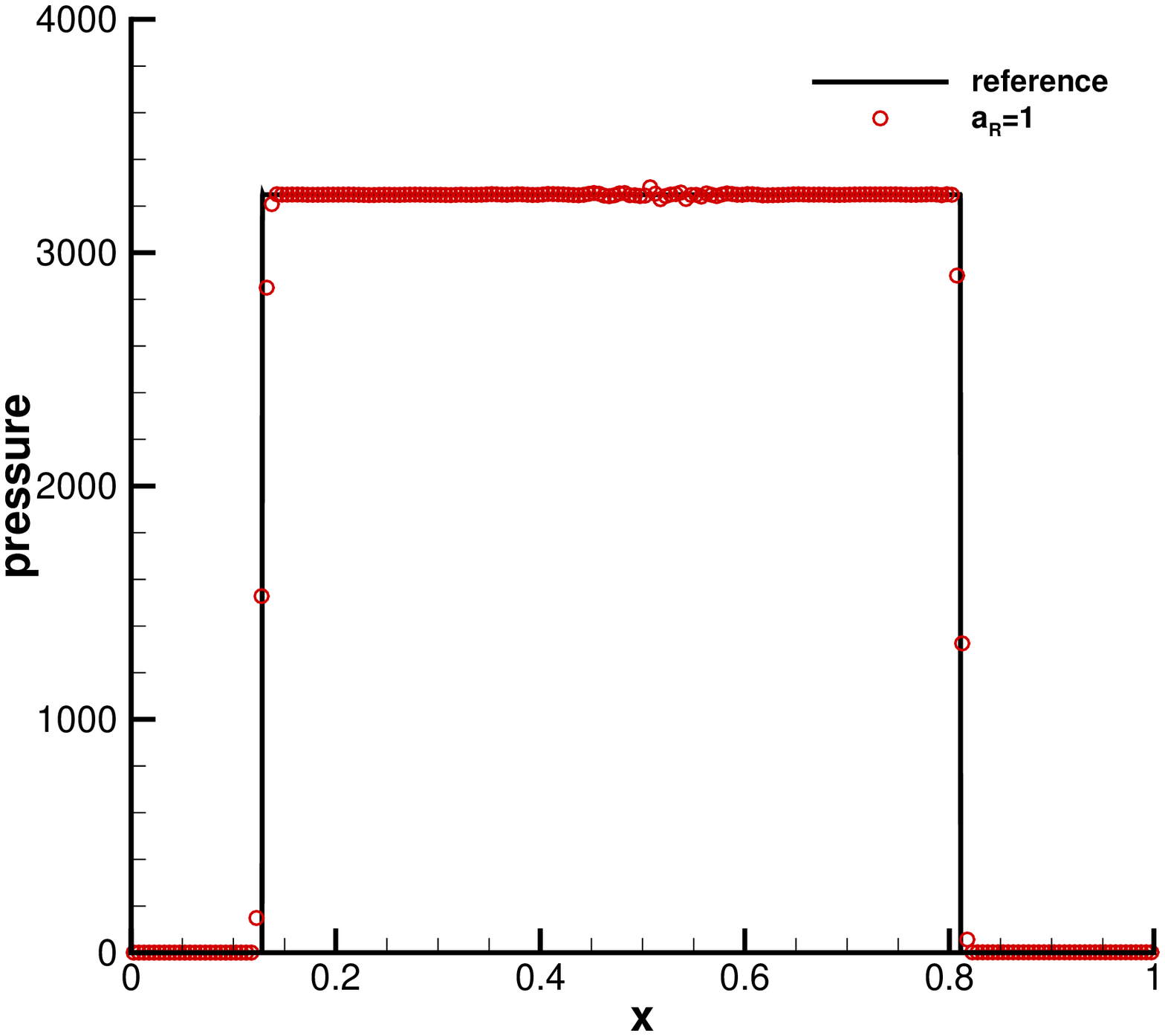}
\caption{\label{riemann-shock-2}  Shock tube problem involving two
shock waves: density, velocity, temperature and pressure at
$t=0.045$ with $a_R=1$ for the numerical results and reference data.
}
\end{figure}

\subsection{Shock tube problem with two rarefaction waves}
In this case, this one-dimensional shock tube problem is tested,
which involves two rarefaction waves moving towards the opposite
directions. The computational domain is $[-0.5,0.5]$, and the
following initial condition is considered
\begin{equation*}
(\rho,U,T) = \begin{cases}
(1,-1, 1),~~  0\leq x<0.5,\\
(1, 1, 1),~~  0.5\leq x\leq1,
\end{cases}
\end{equation*}
where $c_v=1$ and $a_R=1$. For the system without the diffusion
term, i.e., $\kappa=0$, the uniform mesh with 200 cells is used. The
numerical results and the reference solutions are given in
Fig.\ref{riemann-rarefraction-a}, where the reference solutions are
given by the second-order code with 2000 cells. The numerical
results also agree well with the reference solutions
\cite{Radiation-Diffusion-4}. To test the effect of the radiative
diffusion,  the cases with $\kappa=1$, $0.01$ and $0.1$ are tested
as well. The numerical results are shown in
Fig.\ref{riemann-rarefraction-b}, and the numerical results agree
well with \cite{Radiation-Diffusion-4}. With the increase of
$\kappa$, the results reveal that the diffusion effect becomes more
obvious, and such phenomena is quite reasonable in physics.

\begin{figure}[!h]
\centering
\includegraphics[width=0.485\textwidth]{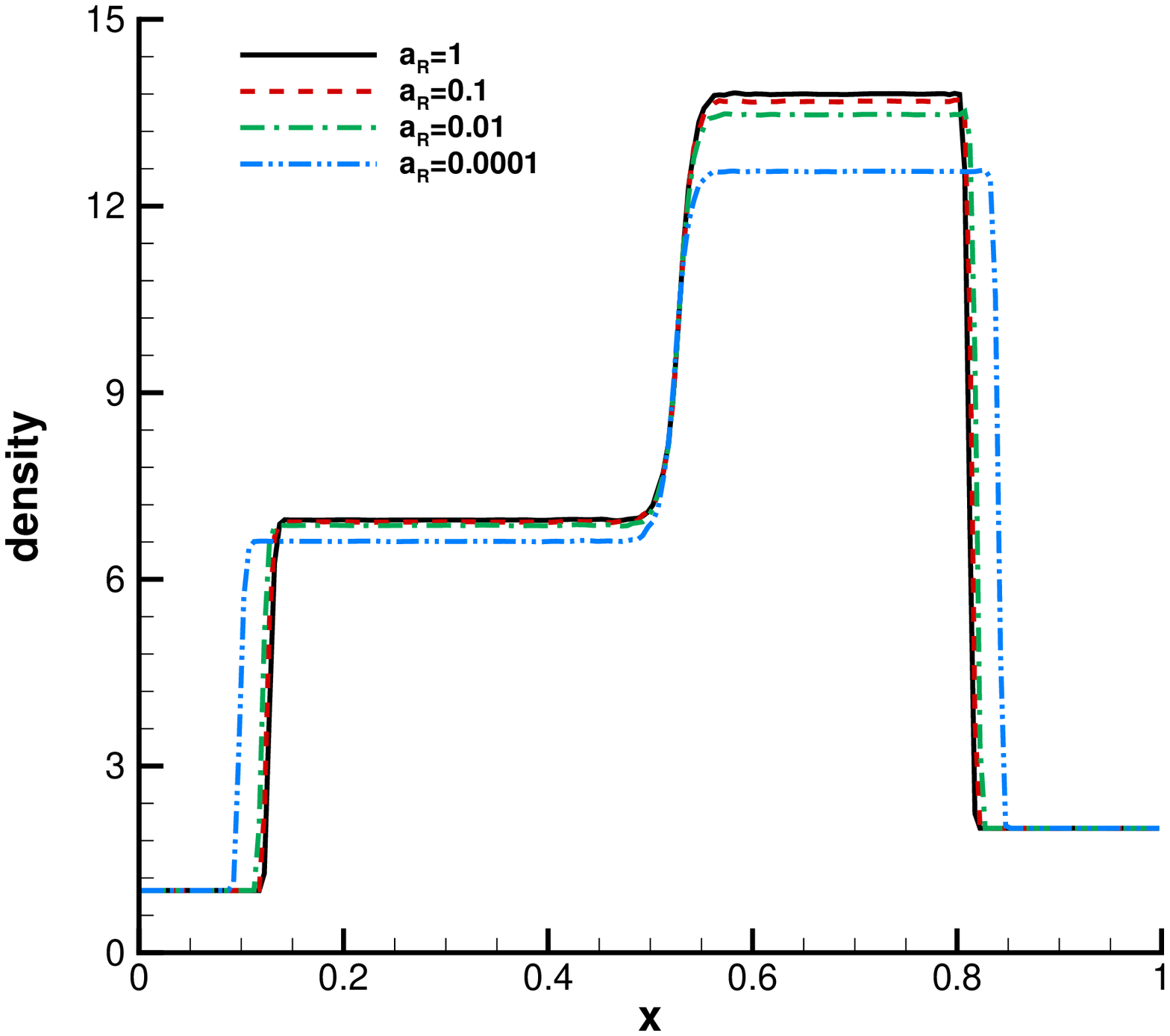}
\includegraphics[width=0.485\textwidth]{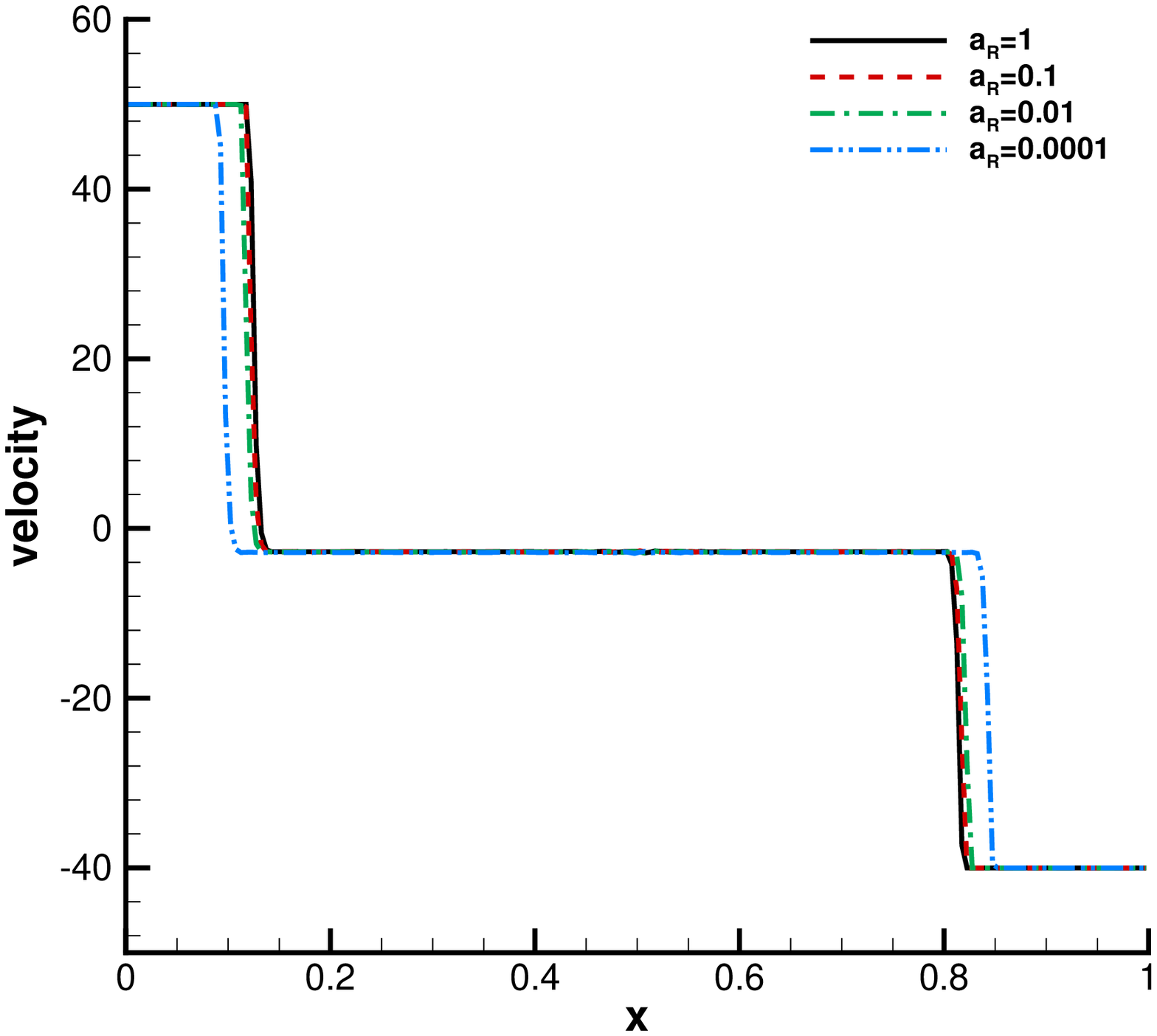}\\
\includegraphics[width=0.485\textwidth]{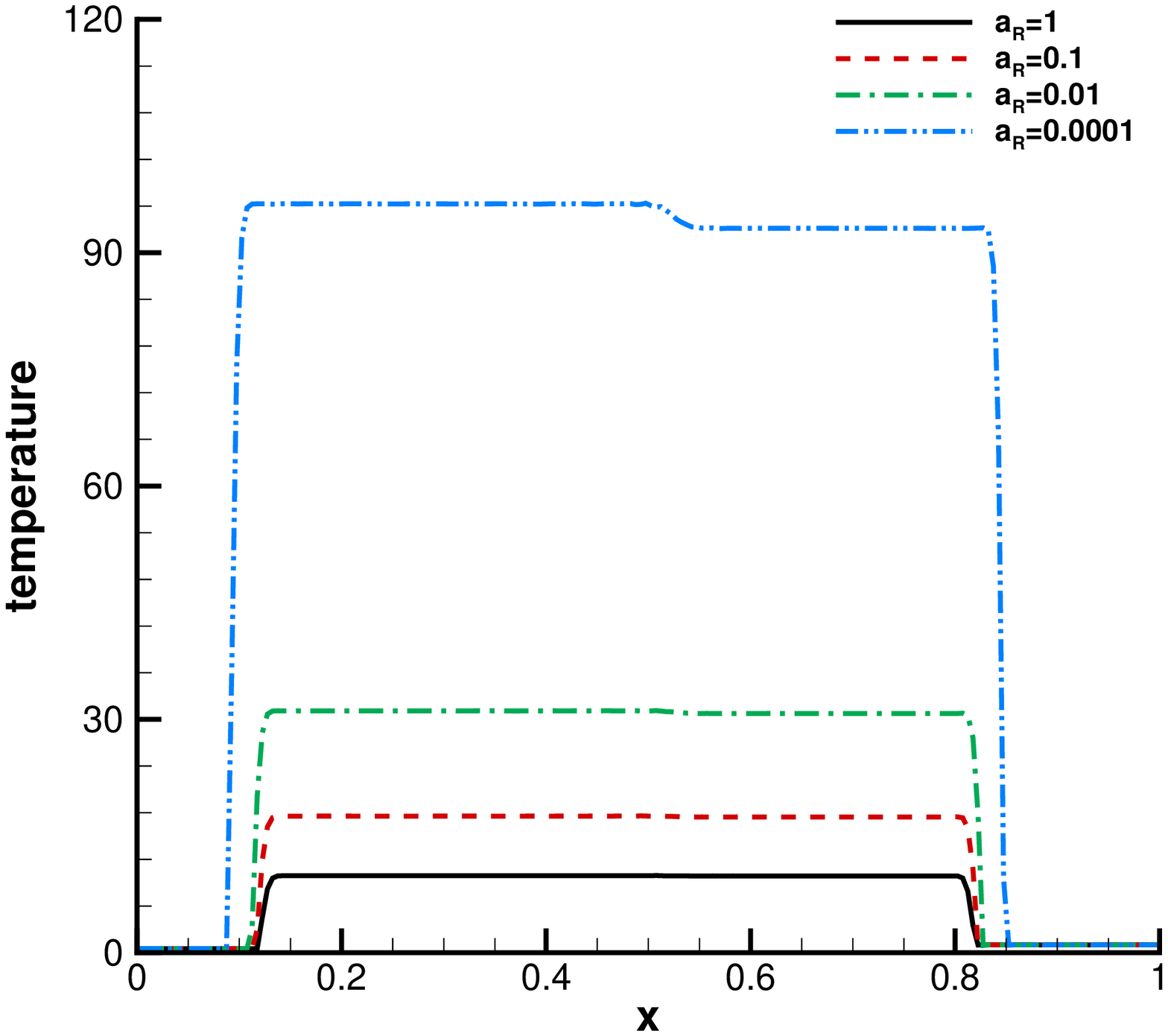}
\includegraphics[width=0.485\textwidth]{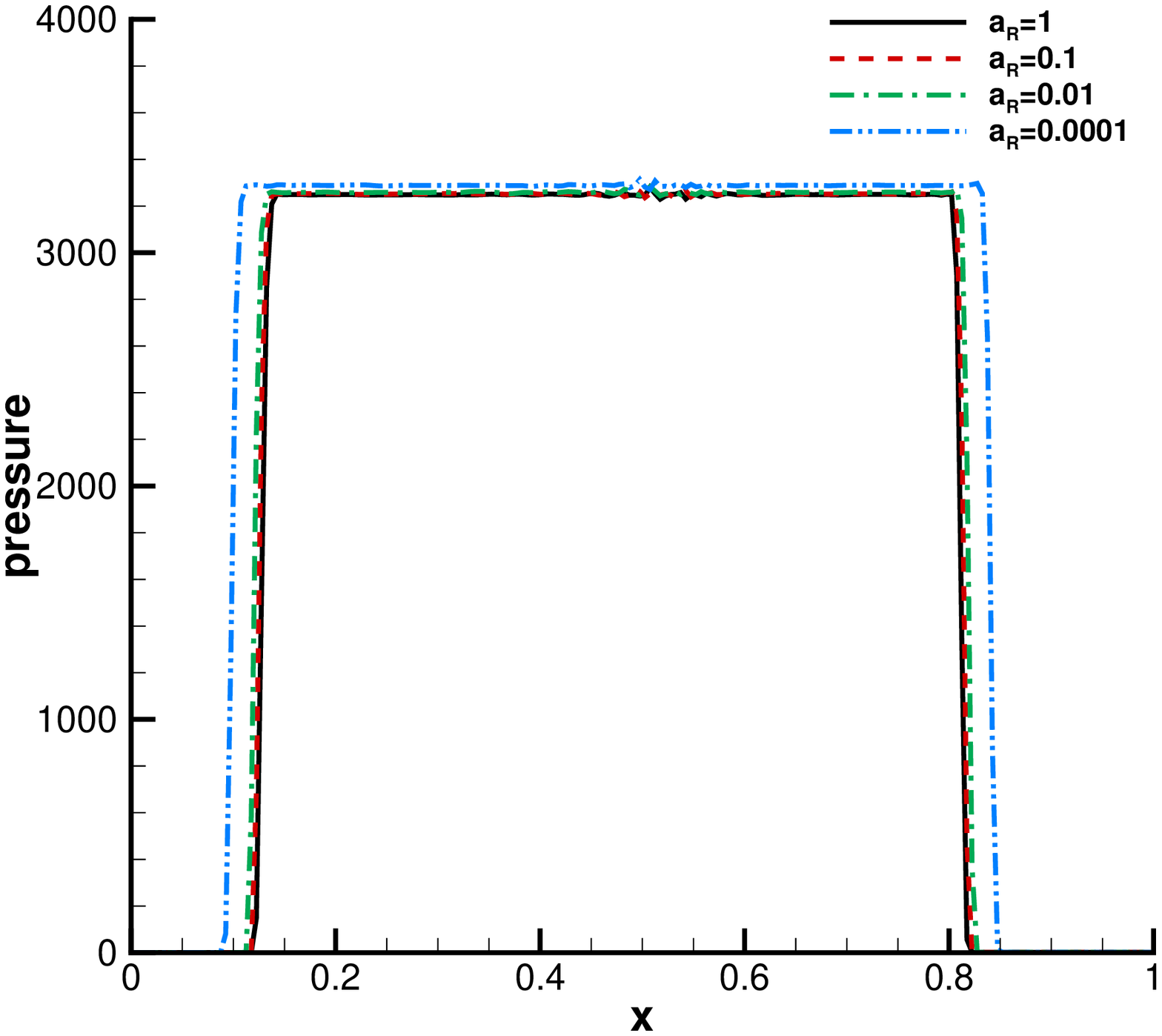}
\caption{\label{riemann-shock-1}  Shock tube problem involving two
shock waves: density, velocity, temperature and pressure at
$t=0.045$ with $a_R=0.0001, 0.01, 0.1$ and $1$.}
\end{figure}

\begin{figure}[!h]
\centering
\includegraphics[width=0.7\textwidth]{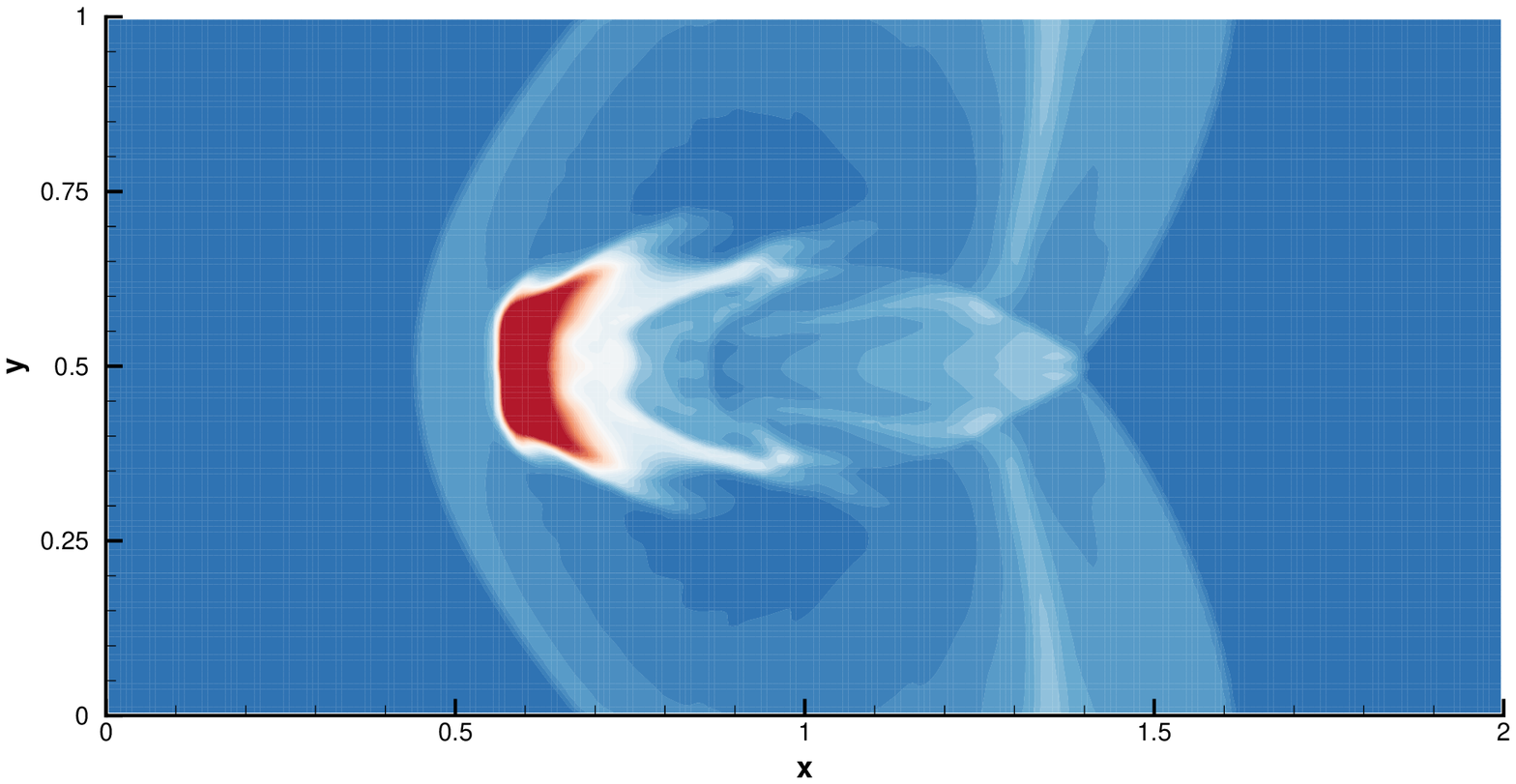}\\
\includegraphics[width=0.7\textwidth]{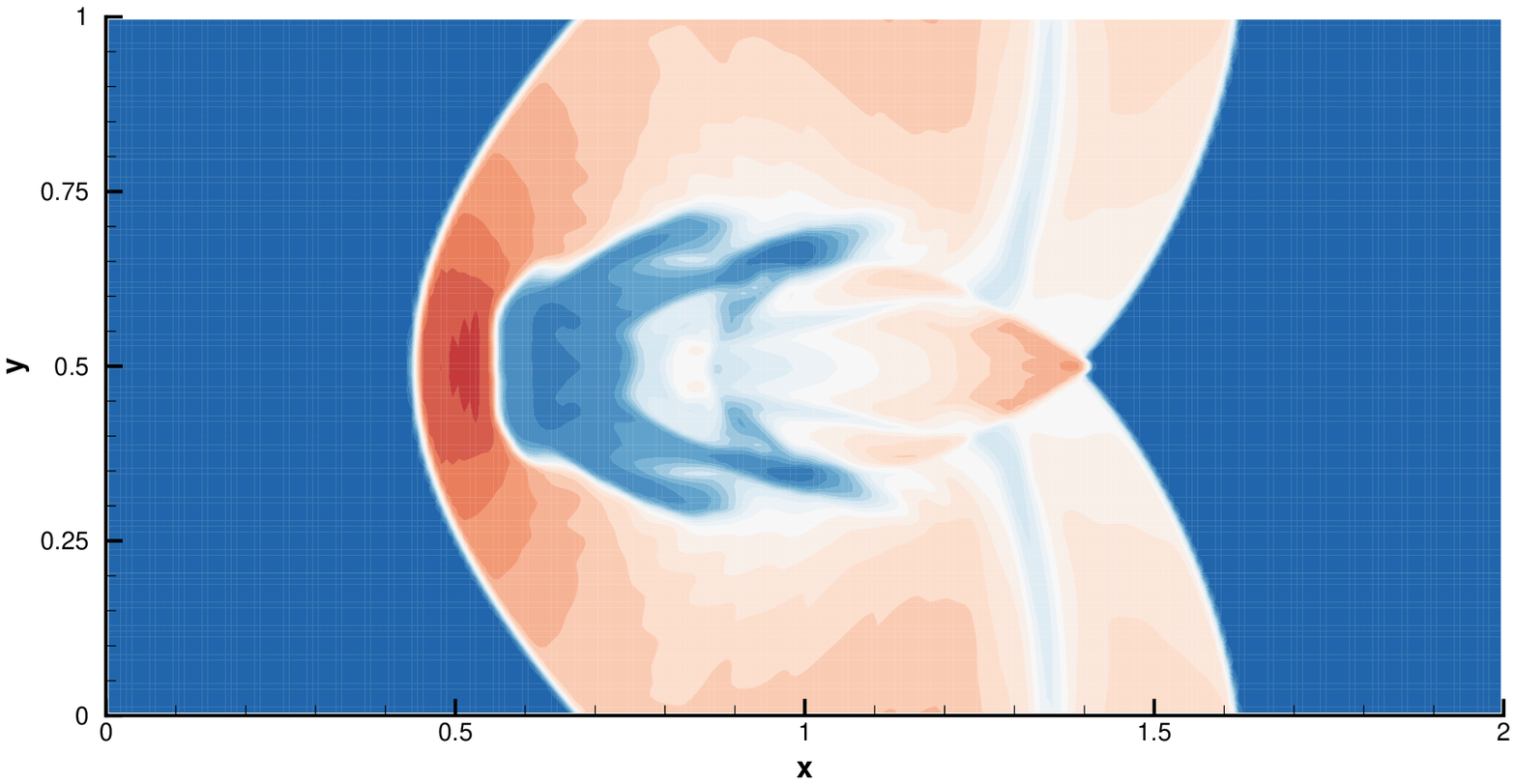}\\
\includegraphics[width=0.7\textwidth]{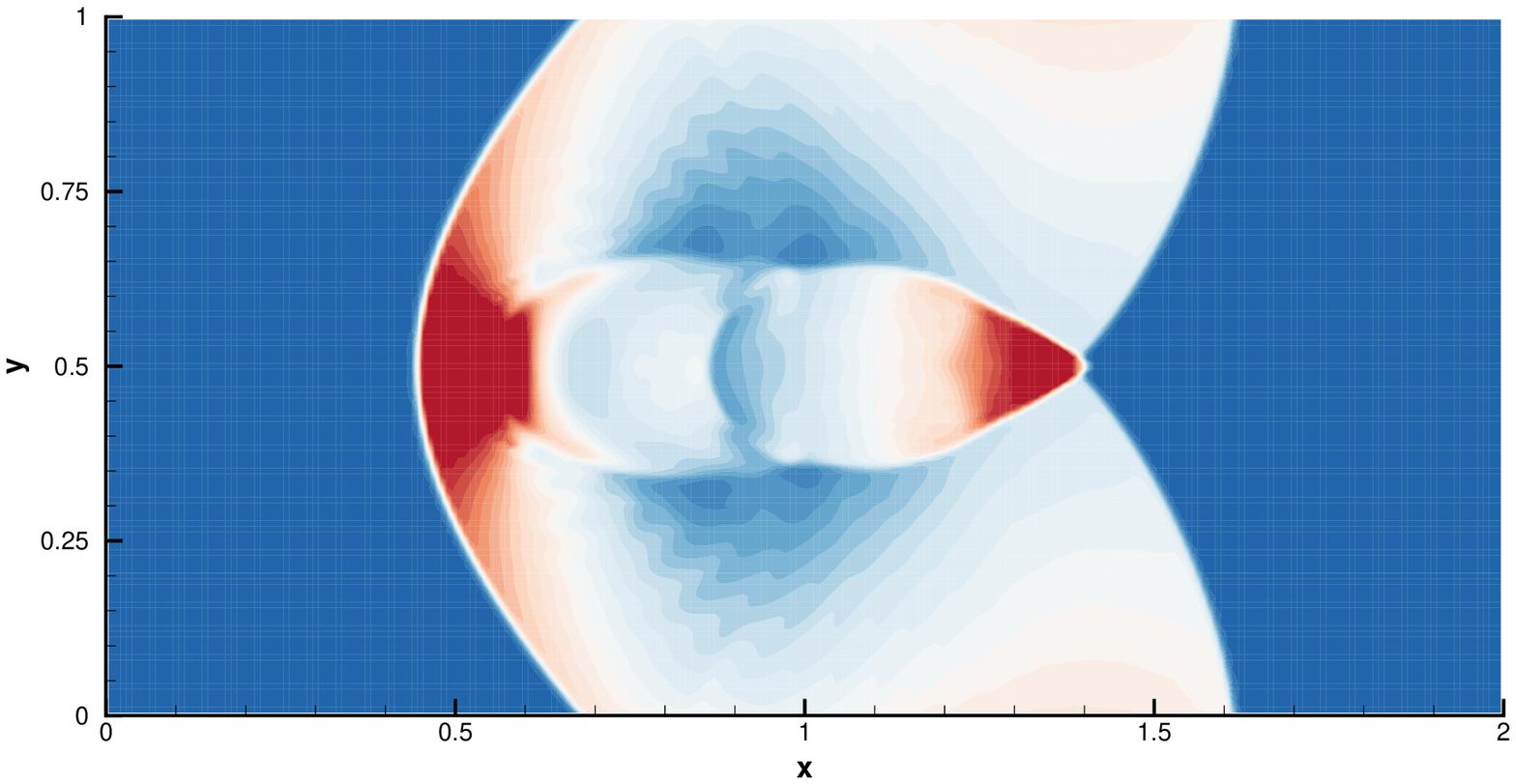}
\caption{\label{bubble0}  Wind and cylindrical bubble interaction: the density,
pressure and temperature distributions at $t = 0.6$ from top to bottom with
with $\kappa(T) = 0$.}
\end{figure}

\subsection{Shock tube problem with two shock waves}
In this case, the one-dimensional shock tube problem involving two
strong shock waves is tested. The computational domain is $[0, 1]$,
and the following initial condition is considered
\begin{equation*}
(\rho,U,T) = \begin{cases}
(1, 50, 0.5),~~   0\leq x\leq 0.65,\\
(2, -40, 1),~~   0.65\leq x\leq 1,
\end{cases}
\end{equation*}
where $c_v=1$, $a_R=1$ and $\kappa=0$. Two shocks with Mach numbers
about $82$ and $39$ are generated from the initial discontinuity.
The uniform mesh with 200 cells are used. The numerical results and
reference solutions at $t=0.045$ are given in
Fig.\ref{riemann-shock-1}, where the reference solutions are given
by the second-order code with 10000 cells. For this case, the exact
solutions for $(\rho, U, T)$ at two sides of contact discontinuity
are $(6.95456, 9.90081, -2.73959)$ and $(13.8008, 9.86595,
-2.73959)$ according to the nonlinear Riemann solver
\cite{Radiation-Diffusion-1}. The numerical solutions are in a good
agreement with the exact solutions. The radiation effect $a_R$ is
tested as well, where $a_R=0.0001, 0.01, 0.1$ and $1$ are used, and
the numerical results are shown in Fig.\ref{riemann-shock-2}. Due to
the strong discontinuity, there is slight oscillation across the
contact discontinuity even with the reconstruction for
characteristic variables \cite{Radiation-Diffusion-3}.

\begin{figure}[!h]
\centering
\includegraphics[width=0.7\textwidth]{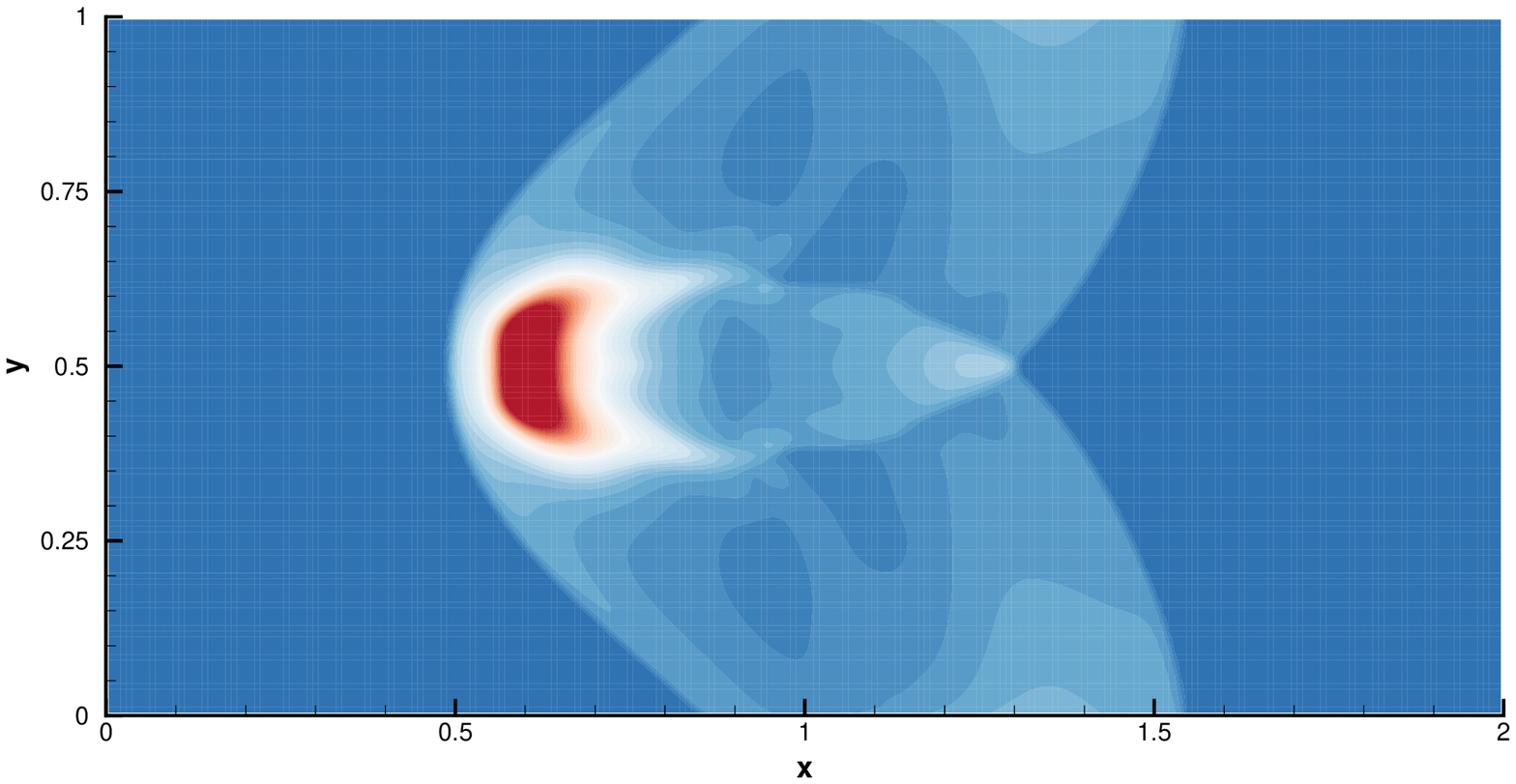}\\
\includegraphics[width=0.7\textwidth]{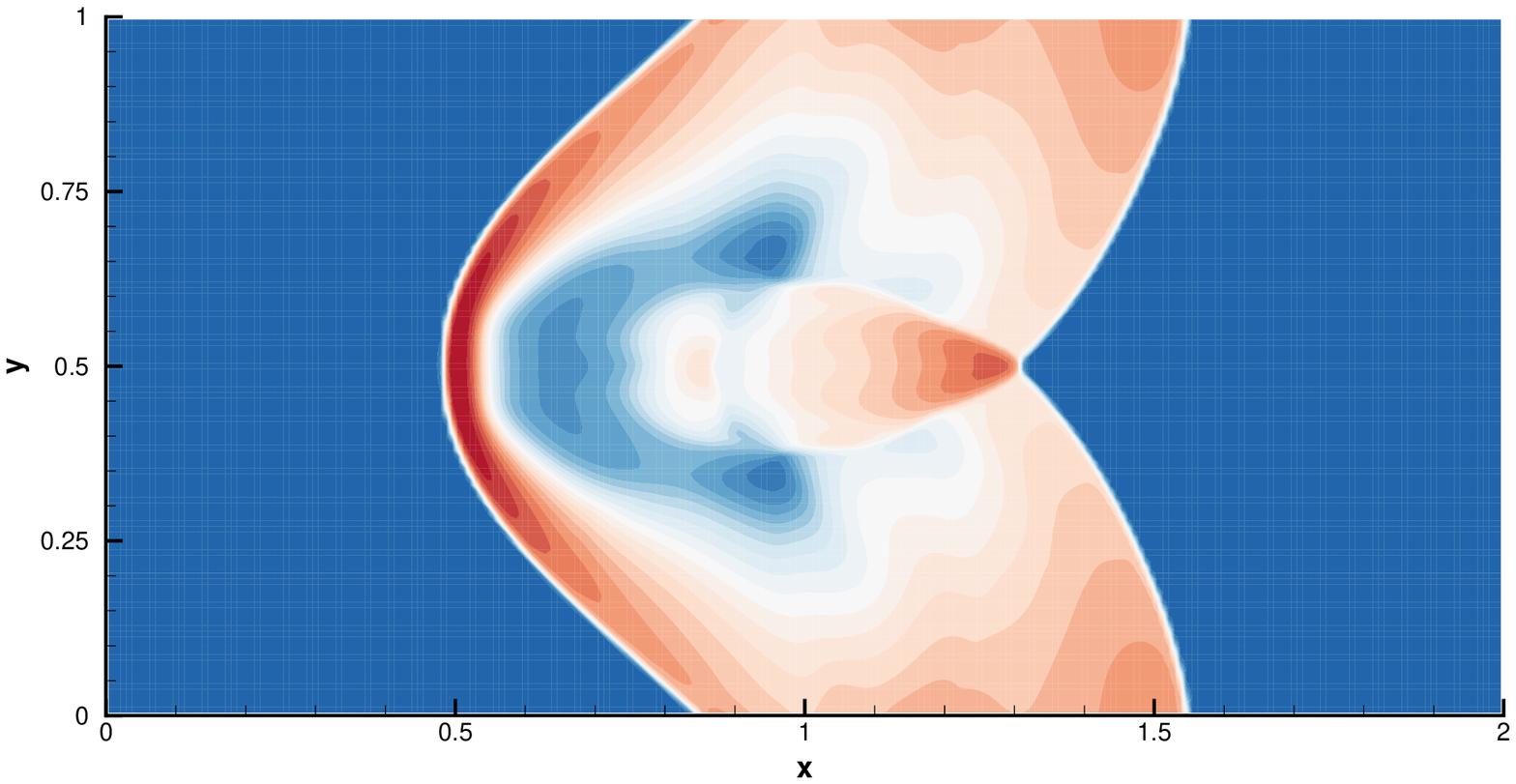}\\
\includegraphics[width=0.7\textwidth]{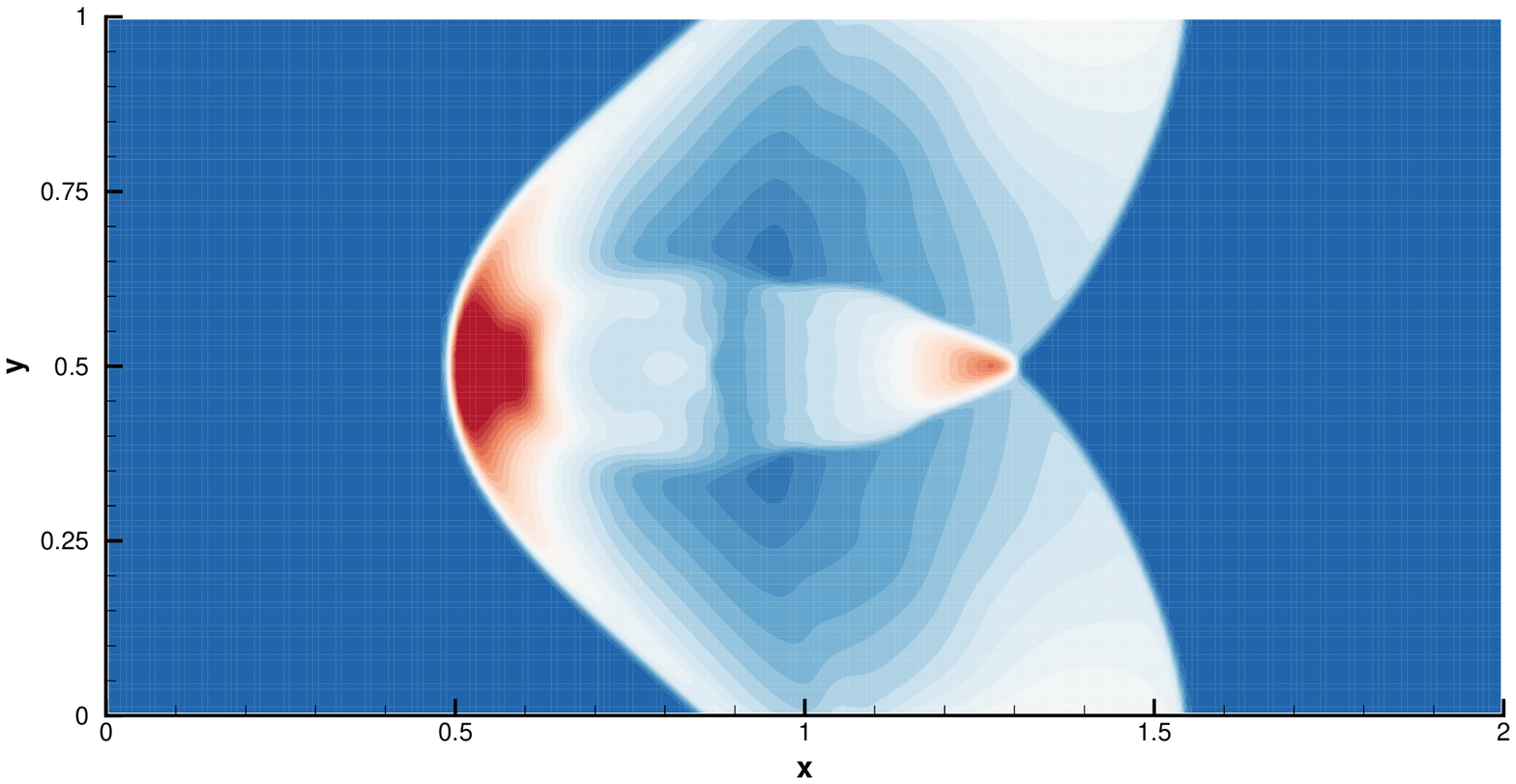}
\caption{\label{bubble1}  Wind and cylindrical bubble interaction: the density,
pressure and temperature distributions at $t = 0.6$ from top to bottom with
with $\kappa(T) = 10^{-3}(1+10T^3)$.}
\end{figure}

\subsection{Interaction between wind and cylindrical bubble}
In this case, the two-dimensional interactions between wind and
denser cylindrical bubble are tested. The simulation is performed in
the domain $[0,2] \times[0,1]$, and there is a cylindrical bubble of
radius $R=0.15$ with its center located at $(0.3,0.5)$. The bubble
is $25$ times denser than the ambient gas, and the temperature of
the cloud is such that the cloud and ambient gas are in an
equilibrium state. Initially, the state for the ambient gas is
$(\rho, U, V, T)=(1, 0, 0, 0.09)$, and the wind is introduced
through the left boundary and assigned as
\begin{equation*}
(\rho, U, V, T)=\big(1, 6(1-e^{-10t}), 0, 0.09\big).
\end{equation*}
The zero gradient boundary condition is given at  the right, upper
and lower boundaries and $a_R$ is taken to be $1$. The computation
is performed with $240\times120$ uniform mesh. For the case with
$\kappa=0$,  the density, pressure and temperature distributions at
$t = 0.6$ are shown in Fig.\ref{bubble0}, where the simulated
results reproduce the large-scale structure of the reference
results, and the high-order scheme resolves the flow structures
better than the second order scheme
\cite{radiative-GKS-0,radiative-GKS-1}. The case with heat
diffusivity $\kappa(T) = 10^{-3}(1+10T^3)$ is also tested. The
density, pressure and temperature distributions at $t = 0.6$ are
shown in Fig.\ref{bubble1}. Due to the heat diffusivity, the flow
structures are smeared.

\subsection{Computational efficiency}
In this case, the computational efficiency of shock tube problem and
interaction between wind and bubble are tested for both
one-dimensional and two-dimensional computation. To show the
efficiency of GMRES procedure,  the case with and without heat
diffusivity are tested, respectively.  In the computation, the
dimension of Krylov subspace in restart procedure is 10, and the
convergence tolerance of GMRES method is  $10^{-8}$.  The CPU times
for different cases are  shown in Table.\ref{tab}  with Intel Core
i7-9700 CPU @ 3.00 GHz within $400$ time steps. The comparison of
CPU time shows the efficiency of implicit-explicit high-order
gas-kinetic scheme with GRMES procedure.

\begin{table}[!h]
\centering
\begin{tabular}{c|c|c}
\hline
~    &  $\kappa$ &CPU time \\
\hline
1D case & $\kappa=0$      &  2.725 \\
1D case & $\kappa=0.01$   &  3.507 \\
2D case & $\kappa=0$      &  81.063 \\
2D case & $\kappa= 10^{-3}(1+10T^3)$  & 266.956  \\
\hline
\end{tabular}
\caption{\label{tab} Computational  efficiency: CPU time (seconds) comparison.}
\end{table}

\section{Conclusion}
In this paper,  a high-order gas-kinetic scheme is proposed for the
equation of radiation hydrodynamic  in the equilibrium-diffusion
limit. Based on the zeroth-order Chapman-Enskog expansion, the
hydrodynamic part of radiation hydrodynamic equation can be obtained
from the modified BGK equation with modified equilibrium state. The
numerical scheme is developed in the finite volume framework, and
the classical multidimensional WENO reconstruction is used to
achieve the spatial accuracy. To achieve the temporal accuracy, a
two-stage method is used,  which is an extension of two-stage
fourth-order method for hyperbolic system. The time scales of
radiation diffusion and hydrodynamic part are different and it will
make the time step of an explicit scheme very small, and an
IMEX-type scheme is developed. The hydrodynamic part is treated
explicitly, and the gas-kinetic solver with the modified equilibrium
state is constructed in the finite volume framework. The nonlinear
Newton-GMRES method is used to treat the radiation diffusion
implicitly,  in which the formation of Jacobian is not required.
One-dimensional and two-dimensional numerical experiments are
carried out, and the numerical results validate the performance of
current scheme.

\section*{Acknowledgements}
The current research of L. Pan is supported by National Natural
Science Foundation of China (11701038) and the Fundamental Research
Funds for the Central Universities, and W.J. Sun is supported by
CAEP foundation (CX20200026) and National Natural Science Foundation
of China (11671048).

\section*{Appendix}
The detailed formulation of matrix $M$ in Eq.\eqref{linear} can be written as
\begin{align*}
M=\left(
\begin{array}{cccccc}
1 & U & V & B'_1&  \displaystyle \frac{K_1}{4 \lambda_1} & \displaystyle\frac{K_2}{4 \lambda_2}  \\
U & U^2+\displaystyle\frac{1}{2\lambda_1} +\displaystyle\frac{1}{2\lambda_2} & UV & B_2'  & \displaystyle \frac{K_1}{4 \lambda_1}U & \displaystyle\frac{K_2}{4 \lambda_2}U  \\
V & UV & V^2+\displaystyle\frac{1}{2\lambda_1} +\displaystyle\frac{1}{2\lambda_2}  & B_3'  & \displaystyle\frac{K_1}{4 \lambda_1}V & \displaystyle\frac{K_2}{4 \lambda_2}V   \\
\displaystyle        B_1 & B_2 & B_3 & B_4  & B_5& B_6\\
\end{array}
\right),
\end{align*}
where
\begin{align*}
    B_1&=\displaystyle\frac{1}{2} (U^2+V^2+\frac{K_1+2}{2\lambda_1}+\frac{K_2+2}{2\lambda_2}),~~
    B_1'=\displaystyle\frac{1}{2} (U^2+V^2+\frac{1}{\lambda_1}+\frac{1}{\lambda_2}),\\
    B_2&=\displaystyle\frac{1}{2} U(U^2+V^2+\frac{K_1+4}{2\lambda_1}+\frac{K_2+4}{2\lambda_2}),~~
    B_2'=\displaystyle\frac{1}{2} U(U^2+V^2+\frac{2}{\lambda_1}+\frac{2}{\lambda_2}),\\
    B_3&=\displaystyle\frac{1}{2}   V(U^2+V^2+\frac{K_1+4}{2\lambda_1}+\frac{K_2+4}{2\lambda_2}),~~
    B_3'=\displaystyle\frac{1}{2} V(U^2+V^2+\frac{2}{\lambda_1}+\frac{2}{\lambda_2}),\\
    B_4&=\displaystyle\frac{1}{4} \big( (U^2+V^2+\frac{1}{\lambda_1}+\frac{1}{\lambda_2})(U^2+V^2+\frac{K_1+2}{2\lambda_1}+\frac{K_2+2}{2\lambda_2})\\
    &+(\frac{2}{\lambda_1}+\frac{2}{\lambda_2})(U^2+V^2+\frac{1}{2\lambda_1}+\frac{1}{2\lambda_2})\big),\\
    B_5&=\displaystyle\frac{1}{4} \big(\frac{K_1}{2\lambda_1}(U^2+V^2)+\frac{K_1^2+4K_1}{4\lambda_1^2}+\frac{K_1(K_2+2)}{4\lambda_1\lambda_2}\big),\\
    B_6&=\displaystyle\frac{1}{4} \big(\frac{K_2}{2\lambda_2}(U^2+V^2)+\frac{K_2^2+4K_2}{4\lambda_2^2}+\frac{K_2(K_1+2)}{4\lambda_1\lambda_2}\big).
\end{align*}


\begin{thebibliography}{10}
\bibitem{Radiation-Diffusion-2}   J.W. Bates, D.A. Knoll, W.J. Rider, R.B. Lowrie, V.A. Mousseauy, On consistent time-integration methods for radiation hydrodynamics in the equilibrium diffusion limit:
Low-energy-density regime, J. Comput. Phys. 167 (2001) 99-130.

\bibitem{BGK-1} P.L. Bhatnagar, E.P. Gross, M. Krook, A Model for Collision Processes in Gases I: Small Amplitude Processes in Charged and Neutral One-Component Systems, Phys. Rev. 94 (1954) 511-525.

\bibitem{Radiation-Diffusion-3} S. Bolding, J. Hansel, J.D. Edwards, J.E. Morel, R.B. Lowrie, Second-order discretization in space and time for radiation-hydrodynamics, J. Comput. Phys. 338 (2017) 511-526.

\bibitem{WENO-Z} R. Borges, M. Carmona, B. Costa, W. S. Don, An improved weighted essentially non-oscillatory scheme for hyperbolic conservation laws, J. Comput. Phys. 227 (2008) 3191-3211.

\bibitem{GMRES-2} P. N. Brown, Y. Saad, Hybrid Krylov Methods for Nonlinear Systems of Equations, SIAM J. Sci. Comput. 11 (1990) 450-481.

\bibitem{Radiation-2} J.I. Castor, Radiation Hydrodynamics, Cambridge University Press (2004).

\bibitem{BGK-2} S. Chapman, T.G. Cowling, The Mathematical theory of Non-Uniform Gases, third edition, Cambridge University Press, (1990).

\bibitem{Radiation-Diffusion-4}  J. Cheng, C.W. Shu, P. Song,  High order conservative Lagrangian schemes for one-dimensional radiation hydrodynamics equations in the equilibrium-diffusion limit, J. Comput.
Phys.  421 (2020) 109724.

\bibitem{Radiation-Diffusion-1}  W. Dai, P.R. Woodward, Numerical simulations for radiation hydrodynamics. I. Diffusion limit, J. Comput. Phys. 142 (1998) 182-207.

\bibitem{GRP-high-2} Z.F. Du, J.Q. Li, A Hermite WENO reconstruction for fourth order temporal accurate schemes based on the GRP solver for hyperbolic conservation laws, J. Comput. Phys. 355 (2018)
385-396.

\bibitem{ENO}  A. Harten, B. Engquist, S. Osher, S.R. Chakravarthy, Uniformly high order accurate essentially non-oscillatory schemes, III, J. Comput. Phys. 71 (1987) 231-303.

\bibitem{UGKS-Xu2}  J. Huang, K. Xu, P. Yu, A unified gas-kinetic scheme for continuum and rarefied flows II: multi-dimensional cases, Commun. Comput. Phys. 12 (2012) 662-690.

\bibitem{GKS-high-2} X. Ji, L. Pan, W. Shyy, K. Xu, A compact fourth-order gas-kinetic scheme for the Euler and Navier-Stokes equations, J. Comput. Phys. 372 (2018) 446-472.

\bibitem{WENO-JS} G.S. Jiang, C.W. Shu, Efficient implementation of weighted ENO schemes, J. Comput. Phys. 126 (1996) 202-228.

\bibitem{radiative-GKS-1} S. Jiang, W.J. Sun, A second-order BGK scheme for the equations of radiation hydrodynamics, Int. J. Numer. Meth. Fluids  53 (2007) 391-416.

\bibitem{Radiation-Diffusion-7}   D.A. Knoll, R.B. Lowrie, J.E Morel, Numerical analysis of time integration errors for non equilibrium radiation diffusion. J. Comput. Phys. 226 (2007) 1332-1347.

\bibitem{GRP-high-1} J.Q. Li, Z.F. Du, A two-stage fourth order time-accurate discretization for Lax-Wendroff type flow solvers I. hyperbolic conservation laws, SIAM J. Sci. Computing, 38 (2016) 3046-3069.

\bibitem{UGKS-Xu3} C. Liu, Y.J. Zhu, K. Xu,  Unified gas-kinetic wave-particle methods I:Continuum and rarefied gas flow, J. Comput. Phys.  401 (2020) 108977.

\bibitem{Radiation-Diffusion-6}   R.B. Lowrie, A comparison of implicit time integration methods for non linear relaxation and diffusion. J. Comput. Phys. 196 (2004) 566-590.

\bibitem{Radiation-Diffusion-5}   R.G. McClarren, T.M. Evans, R.B. Lowrie, J.D. Densmore, Semi-implicit time integration for $P_n$ thermal radiative transfer. J. Comput. Phys. 227 (2008) 7561-7586.

\bibitem{Radiation-1} D. Mihalas, B. W. Mihalas, Foundations of Radiation Hydrodynamics, Oxford University Press (1984).

\bibitem{GKS-high-1} L. Pan, K. Xu, Q.B. Li, J.Q. Li, An efficient and accurate two-stage fourth-order gas-kinetic scheme for the Navier-Stokes equations, J. Comput. Phys.  326 (2016) 197-221.

\bibitem{Radiation-3}  G.C. Pomraning, The Equations of Radiation Hydrodynamics, Pergamon Press, Oxford  (1973).

\bibitem{GMRES-1} Y. Saad, M. H. Schultz, GMRES: A Generalized Minimal Residual Algorithm for Solving Nonsymmetric Linear Systems, SIAM J. Sci. Comput. 7 (1986) 856-869.

\bibitem{radiative-UGKS-1}  W.J. Sun, S. Jiang, K. Xu, An asymptotic preserving unified gas kinetic scheme for gray radiative transfer equations, J. Comput. Phys. 285 (2015) 265-279.

\bibitem{radiative-UGKS-2}  W.J. Sun, S. Jiang, K. Xu, S. Li, An asymptotic preserving unified gas-kinetic scheme for frequency-dependent radiative transfer equations, J. Comput. Phys. 302 (2015) 222-238.

\bibitem{radiative-UGKS-3}  W.J. Sun, S. Jiang, K. Xu,  An implicit unified gas-kinetic scheme for radiative transfer with equilibrium and non-equilibrium diffusive limits, Commun. Comput. Phys. 22 (2015)
899-912.

\bibitem{radiative-UGKS-4}  W.J. Sun, S. Jiang, K. Xu, G.Y. Cao, Multiscale Simulation for the System of Radiation Hydrodynamics, Journal of Scientific Computing (2020) 85:25.

\bibitem{radiative-GKS-2}  W.J. Sun, G.X. Ni, A pressure decoupled BGK model for the equations of radiation hydrodynamics, Journal of Nanjing Normal University, 36 (2013) 5-13.

\bibitem{radiative-GKS-0}  H.Z. Tang, H.M. Wu, Kinetic flux vector splitting for radiation hydrodynamical equations. Computers and Fluids 29 (2000) 917-933.

\bibitem{GKS-Xu1} K. Xu, A gas-kinetic BGK scheme for the Navier-Stokes equations and its connection with artificial dissipation and Godunov method, J. Comput. Phys. 171 (2001) 289-335.

\bibitem{GKS-Xu2} K. Xu, Direct modeling for computational fluid dynamics: construction and application of unified gas kinetic schemes, World Scientific (2015).

\bibitem{UGKS-Xu1}  K. Xu, J. Huang, A unified gas-kinetic scheme for continuum and rarefied flows, J. Comput. Phys. 229 (2010) 7747-7764.
\end{thebibliography}
\end{document}